# Forward-Backward Quantization of Scenario Processes in Multi-Stage Stochastic Optimization

Anna Timonina-Farkas
École Polytechnique Fédérale de Lausanne (EPFL), College of Management of Technology,
Chair of Technology and Operations Management,
EPFL-CDM-TOM, Station 5, CH-1015 Lausanne, Switzerland
International Institute for Management Development (IMD),
Chemin de Bellerive 23, 1003 Lausanne, Switzerland
anna.farkas@epfl.ch, anna.farkas@imd.org

Multi-stage stochastic optimization lies at the core of decision-making under uncertainty. As the analytical solution is available only in exceptional cases, dynamic optimization aims to efficiently find approximations but often neglects non-Markovian time-interdependencies. Methods on scenario trees can represent such interdependencies but are subject to the curse of dimensionality. To ease this problem, researchers typically approximate the uncertainty by smaller but more accurate trees. In this article, we focus on multi-stage optimal tree quantization methods of time-interdependent stochastic processes, for which we develop novel bounds and demonstrate that the upper bound can be minimized via projected gradient descent incorporating the tree structure as linear constraints. Consequently, we propose an efficient quantization procedure, which improves forward-looking samples using a backward step on the tree. We apply the results to the multi-stage inventory control with time-interdependent demand. For the case with one product, we benchmark the approximation because the problem allows a solution in closed-form. For the multi-dimensional problem, our solution found by optimal discrete approximation demonstrates the importance of holding mitigation inventory in different phases of the product life cycle.

## 1. Introduction

Multi-stage stochastic optimization is a well-known technique for solving multi-period decision-making problems in time-interdependent cases. However, the exact solution of a multi-stage optimization program is quite difficult or even impossible to obtain due to the fact that the solutions are functions in a function space (see [2, 31, 37]). As the decisions of multi-period decision-making problems need to be taken in a multi-stage future environment, but the historical time-series data is, clearly, available only for the past, we approximate stochastic processes representing the future via the use of continuous-state probability distributions. Our goal is to study numerical methods for the solution of multi-stage stochastic optimization problems by the use of scenario approximation techniques that are challenging, important and, very often, irreplaceable methods for such optimization problems.

The multi-stage stochastic optimization program, which is formulated in a continuous form, is approximated by the problem, which is finite and discrete. The distance between the problems determines the approximation error. Previously, such distances were only defined for processes in the same probability space. However, the concept of *nested distributions* has combined both the scenario values and the information structure in a purely distributional way, allowing to measure the approximation error between problems defined on different probability spaces [47, 50, 51]. The *nested distance* is introduced in [50, 51, 52] and turns out to be a generalization of the well-known Kantorovich-Wasserstein distance from a single- to a multi-stage case [28, 50, 51, 52, 68]. Also, various contorted nested distances and other metrics, e.g., the weighted measure specifically accounting for heavy-tail distributions (e.g., [3]), utility-based metrics (e.g., [1]), as well as problem-adjusted measures taking objectives and constraints into account (e.g., [29]) are defined for special situations. We focus on the classical definition of the nested distance as in [50, 51, 52].

The fundamental result of A. Pichler and G. Pflug [51] states that the approximation error between problems can be bounded by the term proportional to the nested distance between corresponding nested distributions. Thus, the minimization of the nested distance between the basic continuous model and a finite tree (e.g., with given maximal number of nodes) would lead to the optimal nested distribution that is able to provide fine approximation of the initial problem. However, the minimization of the nested distance between a continuous-state stochastic process and a finitely valued scenario tree is, unfortunately, an intractable optimization problem. Instead, current research focuses on the introduction of various bounds allowing for efficiency in construction of small but accurate trees (see [2, 10, 29, 50, 51, 65, 73]). Following this stream, our article constructs and minimizes the nested distance *upper bounds*. We discretize continuous-state stochastic processes via a novel forward-backward scenario generation algorithm.

The value of the nested distance is strongly dependent on the scenario tree structure. Its slight changes lead to an increase or to a decrease of the nested distance and its upper bounds. In this work, we influence distributions belonging to later stages of scenario trees by making their structure anticipative (i.e., *clairvoyant*) starting at some time stage. A fully clairvoyant tree with completely anticipative scenarios leads to the well-known lower bound on the nested distance, i.e., to the Kantorovich-Wasserstein distance between corresponding multivariate distributions [50, 51]. Moreover, it appears that making the scenarios only partly anticipative improves this lower bound providing a sequence of nested distances bounding the true distance from below. Furthermore, a particular combination of information available on clairvoyant trees results in a novel upper bound for the nested distance combining information about the past with available information about feasible future paths.

Compared to other existing nested distance upper bounds, the bounds proposed in this paper account for both the past scenario realizations and the set of feasible future scenarios at any tree node. Oppositely, the well-known result for the stage-wise optimal quantizers obtained via the minimization of the Kantorovich-Wasserstein distance between distributions at each time stage accounts only for the past scenario realizations and neglects information about feasible future paths [50, 51].

In this article, we introduce probabilistic information about feasible future subtrees explicitly to the scenario generation process. By this, we make a step towards new methods for distributional quantization that are more suitable for multi-stage stochastic optimization programs than existing ones (e.g., stage-wise optimal quantization, Monte-Carlo generation etc.). This is achieved via accounting for the scenario tree structure using linear constraints in the proposed upper bound for the nested distance. Incorporating interdependencies between time stages, we demonstrate that the upper bound can be minimized via a stochastic version of the projected gradient descent [6, 42, 54, 58]. Thus, we develop an efficient quantization algorithm and numerically demonstrate that the resulting quantizers provide better guarantees for the minimal nested distance than the quantizers obtained by the stage-wise minimization of the Kantorovich-Wasserstein distance or, moreover, than the random (Monte-Carlo) sampling.

We use Gaussian and lognormal distributions as baseline models, because these distributions have an explicit conditional form and there are accurate methods that estimate their multivariate probability integrals [62]. Nevertheless, our methods are generalizable to any distributions, of which the conditional form can be estimated based on the historical data. Note that any general distribution can be approximated by a finite Gaussian mixture in a semi-parametric way [18].

We refer to our quantization method as the *forward-backward scenario tree quantization algorithm* and demonstrate its effectiveness for the multi-stage stochastic optimization. We show that it combines high approximation quality with computational efficiency. This is due to the fact that optimal methods allow for using small-size scenario trees in comparison to the necessity of sampling a large number of points in case of Monte-Carlo simulations. Our method takes both the information about the past and the available information about feasible future scenarios into account via a combination of the forward procedure with the backward step on the tree. The backward step starts with the leaves of the tree and goes down to the root taking into account probability distributions on every subtree.

In scenario quantization, the question of computational efficiency is of interest. It is necessary to understand how many values the scenario process should have at each stage of the scenario tree in order to reduce the computational time and to keep the approximation error small. Clearly, one would wish for small scenario trees and still good approximation quality. In this article, we show that even small-size trees can result in a fine approximation of stochastic processes given that the quantization takes the interdependencies between time stages explicitly into account. Our approach is purely distributional and does not depend on the objective function and the constraints of a particular problem (see [29]). Instead, the approximated objective and the constraints depend on the chosen quantizers, guaranteeing the proximity to the initial continuous-state problem. Due to this, our scenario approximation method can be applied in a huge variety of areas: starting with financial planning and inventory control (see [8, 27, 32]), possible applications include topics of energy production, electricity generation planning, supply chain management and similar fields (e.g., [72]). In this article, we focus on multi-stage inventory control with time-interdependent demand.

The multi-stage inventory control problem we consider as an instructive example has a closed-form solution in a one-dimensional case and in a multi-dimensional case with independent product demands. Thus, this gives us an opportunity to compare the approximate solution obtained via our forward-backward algorithm with the true solution. Also, we can assess the quality of the approximation implied by the stage-wise optimal quantization and the Monte-Carlo sampling. Further, we solve the multi-stage inventory control problem in a multi-dimensional case with interdependent products, where we demonstrate the benefits of holding additional inventory across the product life cycle.

The rest of the article is structured as follows: Section 2 describes the mathematical framework, introducing novel bounds on the nested distance. In Section 3, we propose the forward-backward scenario tree quantization algorithm and numerically analyze its functional properties in comparison to the stage-wise optimal quantization and to the Monte-Carlo sampling. Section 4 studies the multi-stage inventory control problem with time-interdependent demand. We consider both cases with a single and multiple product(s).

## 2. Mathematical framework

We start with the probability space $(\Omega, P)$ which carries a filtration $\mathscr{F}_t$, i.e., an increasing sequence of $T$ sigma-algebras $\mathscr{F}_1 \subset \mathscr{F}_2 \subset \ldots, \mathscr{F}_T$. Here $T$ is our decision horizon. On this space, a $D$-dimensional discrete-time continuous-state stochastic process $\xi = (\xi_1, ..., \xi_T)$ is defined, called the *scenario process*, which is adapted to the filtration $\mathscr{F}_t$. We denote this by $\xi \lhd \mathscr{F}_t$. For each individual time step $t$, this means that $x_t$ is measurable with respect to (w.r.t.) $\mathscr{F}_t$, denoted with the same symbol as $x_t \lhd \mathscr{F}_t$. The sigma-algebra $\mathscr{F}_t$ contains sets, which can be observed at stage $t$. Evidently, also our sequence of decisions denoted by $x = (x_0, x_1, \ldots, x_{T-1})$ should be measurable w.r.t the filtration $\mathscr{F}_t$, meaning that $x_t \lhd \mathscr{F}_t$. This statement is referred to as the *non-anticipativity* condition [47, 50, 51]. We call $x$ a *strategy*. The inital decision $x_0$ does not depend on the future observations and is deterministic.

The nested distribution of the scenario process incorporates the values as well as the filtration in a purely distributional manner [47] and allows to measure the distance between such objects without reference to a common probability space. For discretely filtered stochastic processes, the nested distribution is in one-to-one correspondence to a probabilistic tree. However, a filtered stochastic process with continuous-state space also possesses a nested distribution and the distance can - at least conceptually - be defined between continuous-state processes and discrete-state processes, given that they have the same length $T$.

The multi-stage decision problem can be written in the following form [47, 49, 50, 51]:

$$\min_x \left\{ \mathscr{R}\big[Q(x,\xi), \mathscr{F}\big] : x \lhd \mathscr{F}, x \in \mathbb{X} \right\}, \tag{1}$$

where $\mathbb{X}$ is the set of constraints other than the non-anticipativity conditions, $Q(x,\xi)$ is a loss function and $\mathscr{R}(\cdot,\cdot)$ is a functional whose value depends only on the distribution of $Q(x,\xi)$, e.g., an expectation, average-value-at risk etc. [47, 50, 51]. Due to the variational form of the multi-stage stochastic optimization problem (1), its analytical solution can be found explicitly only for the simplest formulations. The solution technique in complex cases is based on the approximation of the problem (1) by a vector optimization problem (2):

$$\min_{\tilde{x}} \left\{ \mathscr{R}\big[Q(\tilde{x},\tilde{\xi}); \widetilde{\mathscr{F}}\big] : \tilde{x} \lhd \widetilde{\mathscr{F}}, \tilde{x} \in \mathbb{X} \right\}, \tag{2}$$

where the stochastic process $\xi$ is replaced by a scenario process with finite number of values $\widetilde{\xi} = (\widetilde{\xi}_1, ..., \widetilde{\xi}_T)$ defined on a finite probability space $(\widetilde{\Omega}, \widetilde{\mathscr{F}}, \widetilde{P})$, which is different from the original infinite one [38, 39, 47].

The distance between problems (1) and (2) determines the approximation error. Denoting nested distributions corresponding to optimization problems (1) and (2) by $\mathbb{P}$ and $\widetilde{\mathbb{P}}$, we employ the notion of the *nested distance* $dl(\mathbb{P}, \widetilde{\mathbb{P}})$, which is introduced in [50, 51]:

DEFINITION 1. (Nested Distance) The nested (multi-stage) distance of order $r \geq 1$ between two nested distributions $\mathbb{P}$ and $\widetilde{\mathbb{P}}$ is defined as

$$
\begin{aligned}
& dl_r(\mathbb{P}, \widetilde{\mathbb{P}}) = \inf_{\pi} \left( \int d(u,v)^r \pi(du,dv) \right)^{\frac{1}{r}}, \qquad (3) \\
& \text{subject to} \quad \mathbb{P} \sim (\Omega, \mathscr{F}, P, \xi),\ \widetilde{\mathbb{P}} \sim (\widetilde{\Omega}, \widetilde{\mathscr{F}}, \widetilde{P}, \widetilde{\xi}), \\
& \pi[A \times \widetilde{\Omega} | \mathscr{F}_t \otimes \widetilde{\mathscr{F}}_t](u,v) = P(A|\mathscr{F}_t)(u),\ (A \in \mathscr{F}_T, 1 \leq t \leq T), \\
& \pi[\Omega \times B | \mathscr{F}_t \otimes \widetilde{\mathscr{F}}_t](u,v) = \widetilde{P}(B|\widetilde{\mathscr{F}}_t)(v),\ (B \in \widetilde{\mathscr{F}}_T, 1 \leq t \leq T),
\end{aligned}
$$

where $d(u,v)$ is the Euclidean distance between transported points $u$ and $v$. Further, we denote by $dl(\mathbb{P}, \widetilde{\mathbb{P}})$ the nested distance of order $r = 1$, i.e. $dl_1(\mathbb{P}, \widetilde{\mathbb{P}}) = dl(\mathbb{P}, \widetilde{\mathbb{P}})$.

In comparison to distance measures defined on the same probability space, the nested distance generalizes the well-known *Kantorovich-Wasserstein distance* from single- to multi-stage problems [28, 50, 51, 68] and allows to account for processes which belong to different probability spaces.

Importantly, in line with the work of G.Ch. Pflug and A.Pichler [50], the nested distributions $\mathbb{P}$ and $\widetilde{\mathbb{P}}$ possess a recursive structure and, hence, the nested distance $dl(\mathbb{P}, \widetilde{\mathbb{P}})$ can be calculated by a recursive procedure. To illustrate this, we let $N_t$ and $\widetilde{N}_t$ be the number of nodes at stage $t$ of scenario trees $\mathbb{P}$ and $\widetilde{\mathbb{P}}$, $\forall t = 1, ..., T$ and suppose that $n_t(u^{t-1}) \in \{1, ..., N_t\}$ and $\widetilde{n}_t(v^{t-1}) \in \{1, ..., \widetilde{N}_t\}$. Note that the number of nodes $N_t$ tends to infinity if the nested distribution $\mathbb{P}$ is continuous. Here, each node $n_t(u^{t-1})$ (resp. $\widetilde{n}_t(v^{t-1})$) depends on the history of the scenario process $u^{t-1}$ (resp. $v^{t-1}$). Further, we let $dl(\mathbb{P}^{t:T}(\cdot|u^{t-1}), \widetilde{\mathbb{P}}^{t:T}(\cdot|v^{t-1}))$ be the nested distance between subtrees $\mathbb{P}^{t:T}(\cdot|u^{t-1})$ and $\widetilde{\mathbb{P}}^{t:T}(\cdot|v^{t-1})$ starting from the nodes $n_t(u^{t-1})$ and $\widetilde{n}_t(v^{t-1})$ at the stage $t-1$. The recursive procedure for the computation (but not the minimization) of the nested distance $dl(\mathbb{P}, \widetilde{\mathbb{P}})$ follows backwards on the scenario trees and can be described via the following two-step procedure:

**Distance between leaves of the tree ($t = T$):** The single-stage nested distance between the final stages of the distributions $\mathbb{P}$ and $\widetilde{\mathbb{P}}$ conditional on the nodes $n_{T-1}(u^{T-1})$ and $\widetilde{n}_{T-1}(v^{T-1})$ is equal to

$$dl\big(\mathbb{P}^{T:T}(\cdot|u^{T-1}), \widetilde{\mathbb{P}}^{T:T}(\cdot|v^{T-1})\big) = d_{KA}\big(P_T(\cdot|u^{T-1}), \widetilde{P}_T(\cdot|v^{T-1})\big),$$

where $d_{KA}\big(P_T(\cdot|u^{T-1}), \widetilde{P}_T(\cdot|v^{T-1})\big)$ is the Kantorovich-Wasserstein distance between the corresponding conditional distributions.

**Nested distance between subtrees ($\forall t = 1,...,T-1$):** At each stage $t = 2,...,T-1$, the nested distance $dl(\mathbb{P}^{t:T}(\cdot|u^{t-1}), \widetilde{\mathbb{P}}^{t:T}(\cdot|v^{t-1}))$ can be computed as the solution of a linear program following from the definition (3) given the distance matrix

$$D_t = \left( d(u_t, v_t) + dl(\mathbb{P}^{t+1:T}(\cdot|u^t), \widetilde{\mathbb{P}}^{t+1:T}(\cdot|v^t)) \right)_{u^t, v^t}, \tag{4}$$

where the distance $dl(\mathbb{P}^{t+1:T}(\cdot|u^t), \widetilde{\mathbb{P}}^{t+1:T}(\cdot|v^t))$ is the outcome of the previous step.

The resulting nested distance $dl(\mathbb{P}, \widetilde{\mathbb{P}})$ is computed at the final iteration with $t = 1$. Importantly, it allows to bound the approximation error between problems (1) and (2) as

$$|V(\mathbb{P}) - V(\widetilde{\mathbb{P}})| \leq C dl(\mathbb{P}, \widetilde{\mathbb{P}}), \tag{5}$$

where $V(\mathbb{P})$ and $V(\widetilde{\mathbb{P}})$ are the respective optimal solutions and $C$ is a constant depending on the properties of the optimization problem. In order to reduce the approximation error, one aims for the minimization of the nested distance $dl(\mathbb{P}, \widetilde{\mathbb{P}})$. However, even the computation of the nested distance $dl(\mathbb{P}, \widetilde{\mathbb{P}})$ is a complex procedure if the scenario tree $\mathbb{P}$ is large or, moreover, continuous. The complexity only increases for the direct minimization of the nested distance over $\widetilde{\mathbb{P}}$. As a solution to this issue, bounds on the nested distance are introduced and the upper bounds are minimized for the near-optimal approximation.

On one hand, the well-known lower bound on the nested distance is implied by the Kantorovich-Wasserstein distance $d_{KA}(P, \widetilde{P})$ between the joint distribution $P$ of the stochastic process $\xi$ and its discrete approximation $\widetilde{P}$. On the other hand, the sum of supremal Kantorovich-Wasserstein distances between the conditional distributions $P_t(\cdot|u^{t-1})$ sitting at stages $1,...,T$ provides the upper bound on the nested distance [50, 51], i.e.,

$$d_{KA}(P, \widetilde{P}) \leq dl(\mathbb{P}, \widetilde{\mathbb{P}}) \leq \sum_{t=1}^{T} \sup_{u^{t-1}, v^{t-1}} d_{KA}\big(P_t(\cdot|u^{t-1}), \widetilde{P}_t(\cdot|v^{t-1})\big). \tag{6}$$

Here, the lower bound relaxes the tree structure and does not account for gradually increasing information. The upper bound depends solely on stage-wise conditional distributions limiting the probabilistic information about all future stages. Clearly, it does not take the recursive structure of the nested distance into account neglecting how far from each other the subtrees $\mathbb{P}^{t:T}(\cdot|u^{t-1})$ and $\widetilde{\mathbb{P}}^{t:T}(\cdot|v^{t-1})$ are.

Moreover, if the Lipschitz property holds for the distribution $P_t(\cdot|u^{t-1})$ with constants $K_t$ $\forall t$, i.e., if $d_{KA}\big(P_t(\cdot|u^{t-1}), P_t(\cdot|v^{t-1})\big) \leq K_t d(u^{t-1}, v^{t-1})$, then

$$dl(\mathbb{P}, \widetilde{\mathbb{P}}) \leq \sum_{t=1}^{T} d_{KA}(P_t, \widetilde{P}_t) \prod_{s=t+1}^{T} (K_t + 1), \tag{7}$$

where $d_{KA}(P_t, \widetilde{P}_t) = \sup_{v^{t-1}} d_{KA}\big(P_t(\cdot|v^{t-1}), \widetilde{P}_t(\cdot|v^{t-1})\big)$. In order to incorporate both (i) the recursive structure of scenario trees and (ii) the probabilistic information of nested distributions, we further state new lower and upper bounds.

### 2.1. Recursive nested distance bounds

Direct minimization of the nested distance $dl(\mathbb{P}, \widetilde{\mathbb{P}})$ would lead to a fine result in the quality of approximation of (1) by (2), as the difference between $V(\mathbb{P})$ and $V(\widetilde{\mathbb{P}})$ would be bounded. However, the minimization of the nested distance between a continuous stochastic process and a scenario tree is a complex task as the transportation plan $\pi$ is a matrix of an infinite size for the case of continuous stochastic processes. Instead of searching for the optimal $\pi$, we propose lower and upper bounds for the nested distance which lead to efficient scenario quantization algorithms. For this, we introduce the concept of a $t$-*clairvoyant* tree, which we denote by $\mathbb{P}_{c(t)}$.

DEFINITION 2. ($t$ - Clairvoyant tree) A scenario tree $\mathbb{P}_{c(t)}$ is called *clairvoyant* at/from stage $t$ if the filtration structure starting with stage $t$ is equal to $\mathscr{F}_T$, i.e., to the filtration at the final stage $T$.

Clearly, every tree $\mathbb{P}$ can be expanded to the clairvoyant form $\mathbb{P}_{c(t)}$ from any stage $t$ by changing its structure. Starting at the stage $t$, all scenarios $\widetilde{\xi}^{t:T}_{c(t)} = (\widetilde{\xi}_{c(t),t}, ..., \widetilde{\xi}_{c(t),T})$ of the tree $\mathbb{P}_{c(t)}$ can be obtained from scenarios $\widetilde{\xi}^{t:T} = (\widetilde{\xi}_t, ..., \widetilde{\xi}_T)$ of the tree $\mathbb{P}$ by relaxing its filtration. This corresponds to the case demonstrated in Figure 1 given $x_1 = x_1' = x_1''$. All random variables $\xi_t, \xi_{t+1}, ..., \xi_T$ are measurable with respect to the filtration $\mathscr{F}_T$.

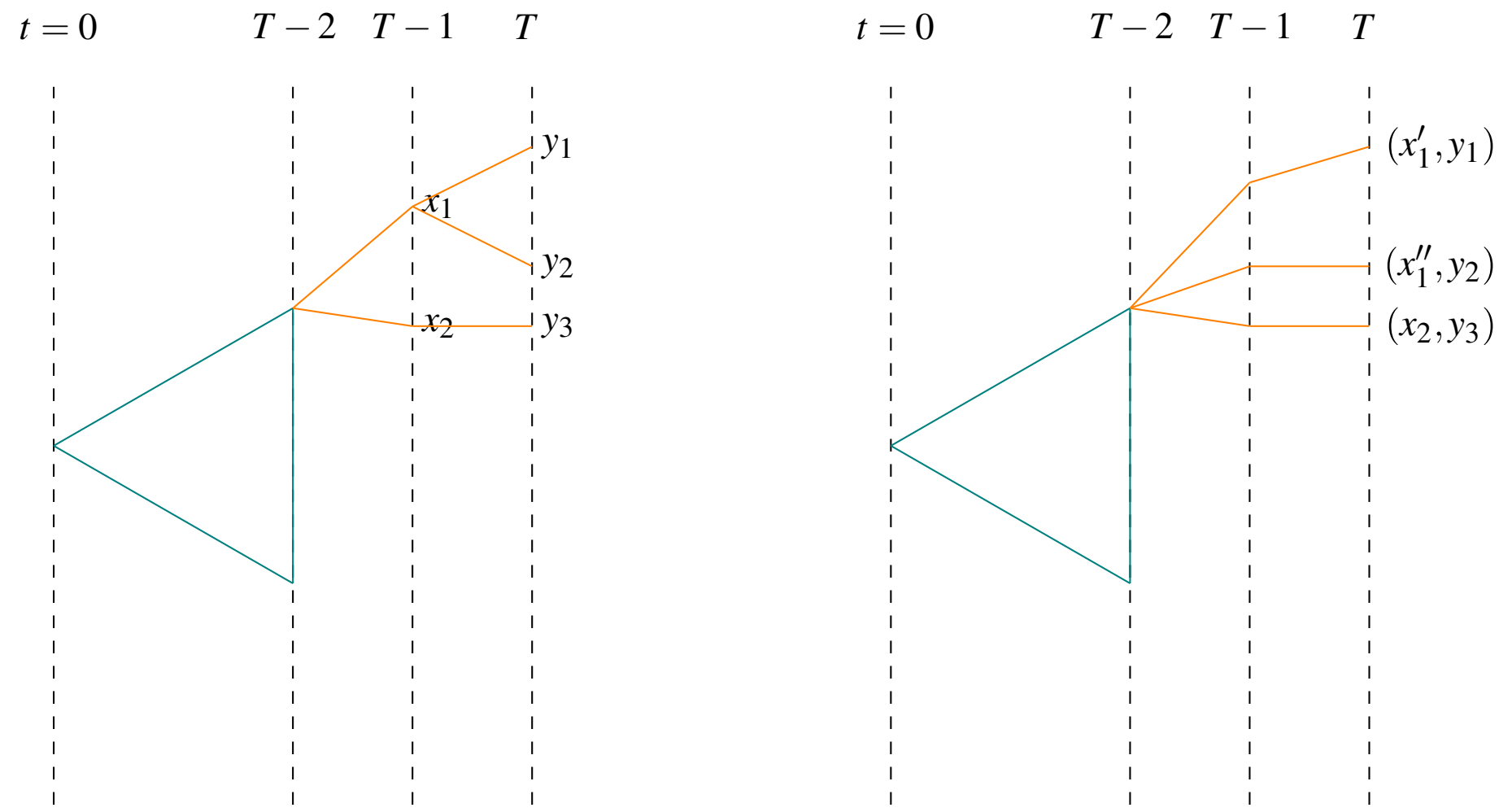

**Figure 1** **Non-clairvoyant v.s. clairvoyant tree at the stage $T-1$.**

Note that some scenario values at stages $t, ..., T-1$ of the clairvoyant tree $\mathbb{P}_{c(t)}$ may coincide with each other, e.g., if $x_1 = x_1' = x_1''$ in Figure 1. In this case, the joint distribution $P^{t:T}$ of random variables at stages $t, ..., T$ of the clairvoyant tree $\mathbb{P}_{c(t)}$ coincides with the corresponding multivariate distribution of the tree $\mathbb{P}$ [50]. We refer to this case as the *clairvoyant representation* of a given tree $\mathbb{P}$. Differently, if some values and corresponding probabilities are different (e.g., if $x_1 \neq x_1' \neq x_1''$ in Figure 1), the distributions are obviously not equal. In both cases, the distance between nested distributions $\mathbb{P}$ and $\mathbb{P}_{c(t)}$ is not zero due to the difference in filtration structure (see Example 2.6 of [22], Example 3.10 of [50]).

LEMMA 1. (Lower Bounds) Let $\mathbb{P}$ and $\widetilde{\mathbb{P}}$ be two nested distributions with corresponding multivariate distributions $P$ and $\widetilde{P}$ and suppose that their *clairvoyant representations* $\mathbb{P}_{c(t)}$ and $\widetilde{\mathbb{P}}_{c(t)}$ are defined respectively $\forall t = 1,...,T$. Then, the following chain of lower bounds holds:

$$d_{KA}(P,\widetilde{P}) = dl(\mathbb{P}_{c(1)},\widetilde{\mathbb{P}}_{c(1)}) \leq ... \leq dl(\mathbb{P}_{c(t)},\widetilde{\mathbb{P}}_{c(t)}) \leq ... \leq dl(\mathbb{P}_{c(T)},\widetilde{\mathbb{P}}_{c(T)}) = dl(\mathbb{P},\widetilde{\mathbb{P}}). \quad (8)$$

**Proof:** Note that the nested distance between $\mathbb{P}_{c(1)}$ and $\widetilde{\mathbb{P}}_{c(1)}$ coincides with the Kantorovich-Wasserstein distance between corresponding multivariate distributions, i.e.,

$$d_{KA}(P,\widetilde{P}) = dl(\mathbb{P}_{c(1)},\widetilde{\mathbb{P}}_{c(1)}).$$

This directly follows from the nested distance definition (see [50, 51]), which is equivalent to the Kantorovich-Wasserstein distance definition in case filtration conditions are relaxed.

Considering all subtrees of $\mathbb{P}$ and $\widetilde{\mathbb{P}}$ at the stage $t$ and relaxing the filtration structure of each of them (i.e., receiving $\mathbb{P}_{c(t)}$ and $\widetilde{\mathbb{P}}_{c(t)}$), we decrease the nested distance between these subtrees to the Kantorovich-Wasserstein distance between corresponding multivariate distributions. Therefore, the resulting nested distance $dl(\mathbb{P},\widetilde{\mathbb{P}})$ also decreases due to its recursive structure [50]. Also, we claim that $dl(\mathbb{P}_{c(t-1)},\widetilde{\mathbb{P}}_{c(t-1)}) \leq dl(\mathbb{P}_{c(t)},\widetilde{\mathbb{P}}_{c(t)})$, taking into account that the filtration relaxation at the stage $t-1$ automatically implies the filtration relaxation at the stage $t$. Combining the inequalities for stages $1,...,T$ leads to Lemma 1. □

Importantly, Lemma 1 holds true only if both nested distributions $\mathbb{P}$ and $\widetilde{\mathbb{P}}$ are decomposed to clairvoyant trees at stages $t = 1,...,T$. Counterexample 2 presented in Appendix 6 shows that making clairvoyant only one of two trees $\mathbb{P}$ or $\widetilde{\mathbb{P}}$ may lead to a decrease or to an increase in the nested distance. Also note that computing nested distances in statement (8) requires a recursive procedure described in [50]: one starts with Kantorovich-Wasserstein distances between conditional distributions sitting at the stage $T$ and constructs a distance matrix for a linear program (LP) to proceed to the stage $T-1$. At every stage $1,...,T-1$, the nested distance can be computed as the solution of a linear program with a distance matrix constructed at the previous stage (see [50]).

Next, we introduce several upper bounds on the nested distance $dl(\mathbb{P},\widetilde{\mathbb{P}})$:

THEOREM 1. (Upper Bound) Let $\mathbb{P}$ and $\widetilde{\mathbb{P}}$ be two nested distributions. The following sequence of upper bounds holds $\forall t = 1,...,T$ and $s = 1,...,t$:

$$dl_t(\mathbb{P},\widetilde{\mathbb{P}}) \leq dl_{t-s}(\mathbb{P},\widetilde{\mathbb{P}}) + \sup_{u^{t-s},v^{t-s}} dl^{t-s+1:t}(\mathbb{P}(\cdot|u^{t-s}),\widetilde{\mathbb{P}}(\cdot|v^{t-s})). \quad (9)$$

**Proof:** Denote the nested distance between two trees that are cut beyond time $t$ by $dl_t(\mathbb{P},\widetilde{\mathbb{P}})$ $\forall t$. The stages $t+1,...,T$ of distributions $\mathbb{P}$ and $\widetilde{\mathbb{P}}$ are not taken into account in the nested distance $dl_t(\mathbb{P},\widetilde{\mathbb{P}})$. According to the definition of the nested distance we can write $dl_t(\mathbb{P},\widetilde{\mathbb{P}}) \leq \int d(u^t,v^t)\pi(du^t,dv^t)$, $\forall \pi \in \Pi \quad \forall t = 1,...,T$, where $d(u^t,v^t) = \sum_{i=1}^{t} d(u_i,v_i) = d(u^{t-s},v^{t-s}) + d(u^{t-s+1:t},v^{t-s+1:t})$, $\forall s = 1,..,t$ and $\Pi$ is the feasible set.

Denote by $\pi^*(u^{t-s}, v^{t-s})$ the transportation plan optimal for the nested distance $dl_{t-s}(\mathbb{P}, \widetilde{\mathbb{P}})$ and by $\pi^*(u^{t-s+1:t}, v^{t-s+1:t}|u^{t-s}, v^{t-s})$ the transportation plan optimal for the nested distance $dl^{t-s+1:t}(\mathbb{P}(\cdot|u^{t-s}), \widetilde{\mathbb{P}}(\cdot|v^{t-s}))$. Using the feasible transportation plan

$$\pi(du^t, dv^t) = \pi^*(u^{t-s+1:t}, v^{t-s+1:t}|u^{t-s}, v^{t-s})\pi^*(du^{t-s}, dv^{t-s}) \in \Pi$$

we bound the nested distance as

$$\begin{aligned} dl_t(\mathbb{P}, \widetilde{\mathbb{P}}) \leq & \int [d(u^{t-s}, v^{t-s}) + d(u^{t-s+1:t}, v^{t-s+1:t})]\pi^*(u^{t-s+1:t}, v^{t-s+1:t}|u^{t-s}, v^{t-s})\pi^*(du^{t-s}, dv^{t-s}) = \\ & = \int d(u^{t-s}, v^{t-s})\pi^*(du^{t-s}, dv^{t-s}) + \int d(u^{t-s+1:t}, v^{t-s+1:t})\pi^*(u^{t-s+1:t}, v^{t-s+1:t}|u^{t-s}, v^{t-s})\pi^*(du^{t-s}, dv^{t-s}) = \\ & = dl_{t-s}(\mathbb{P}, \widetilde{\mathbb{P}}) + \int dl^{t-s+1:t}(\mathbb{P}(\cdot|u^{t-s}), \widetilde{\mathbb{P}}(\cdot|v^{t-s}))\pi^*(du^{t-s}, dv^{t-s}), \end{aligned}$$

where $dl^{t-s+1:t}(\mathbb{P}(\cdot|u^{t-s}), \widetilde{\mathbb{P}}(\cdot|v^{t-s})) \leq \sup_{u^{t-s}, v^{t-s}} dl^{t-s+1:t}(\mathbb{P}(\cdot|u^{t-s}), \widetilde{\mathbb{P}}(\cdot|v^{t-s}))$ and, hence,

$$dl_t(\mathbb{P}, \widetilde{\mathbb{P}}) \leq dl_{t-s}(\mathbb{P}, \widetilde{\mathbb{P}}) + \sup_{u^{t-s}, v^{t-s}} dl^{t-s+1:t}(\mathbb{P}(\cdot|u^{t-s}), \widetilde{\mathbb{P}}(\cdot|v^{t-s})), \ \forall t = 1, ..., T; \ s = 1, ..., t.$$

If $T$ is greater than two and is even, one can subdivide the scenario tree into two-stage parts and write:

$$dl(\mathbb{P}, \widetilde{\mathbb{P}}) \leq \sum_{j=2i,\, i=1}^{i=T/2} \sup_{u^{j-2}, v^{j-2}} dl^{j-1:j}(\mathbb{P}(\cdot|u^{j-2}), \widetilde{\mathbb{P}}(\cdot|v^{j-2})). \tag{10}$$

If the number of tree stages is odd, the bound is

$$dl(\mathbb{P}, \widetilde{\mathbb{P}}) \leq \sum_{j=2i,\, i=1}^{i=(T-1)/2} \sup_{u^{j-2}, v^{j-2}} dl^{j-1:j}(\mathbb{P}(\cdot|u^{j-2}), \widetilde{\mathbb{P}}(\cdot|v^{j-2})) + \sup_{u^{T-1}, v^{T-1}} d_{KA}(P_T(\cdot|u^{T-1}), \widetilde{P}_T(\cdot|v^{T-1})). \tag{11}$$

□

Note that both upper bounds (6) and (7) are an immediate consequence of the bound (9) given $s = t-1$. Further, we can state Corollary (1) using the concept of clairvoyant scenario trees. For this, we use multivariate conditional probabilities $P^{t:T}(\cdot|\mathscr{F}_T(u^{t-1}))$ and $\widetilde{P}^{t:T}(\cdot|\mathscr{F}_T(v^{t-1}))$, which are the joint distributions between stages $t$ and $T$ of scenario trees $\mathbb{P}_{c(t)}$ and $\widetilde{\mathbb{P}}_{c(t)}$ given past scenarios $u^{t-1}$ and $v^{t-1}$ and the filtration structure $\mathscr{F}_T$.

COROLLARY 1. Consider clairvoyant scenario trees $\mathbb{P}_{c(t)}$ and $\widetilde{\mathbb{P}}_{c(t)}$ $\forall t = 1, ..., T$ with multivariate conditional probabilities $P^{t:T}(\cdot|\mathscr{F}_T(u^{t-1}))$ and $\widetilde{P}^{t:T}(\cdot|\mathscr{F}_T(v^{t-1}))$. The following upper bound holds:

$$\inf_{\widetilde{\mathbb{P}}} dl(\mathbb{P}, \widetilde{\mathbb{P}}) \leq \sum_{t=1}^{T} \inf_{\widetilde{P}^{t:T}} \sup_{u^{t-1}, v^{t-1}} d_{KA}\big(P^{t:T}(\cdot|\mathscr{F}_T(u^{t-1})), \widetilde{P}^{t:T}(\cdot|\mathscr{F}_T(v^{t-1}))\big). \tag{12}$$

**Proof:** Using the definition of the nested distance we can write $dl_t(\mathbb{P}, \widetilde{\mathbb{P}}) \leq \int d(u^t, v^t)\pi(du^t, dv^t)$, $\forall \pi \in \Pi \quad \forall t = 1, ..., T$, where $d(u^t, v^t) = \sum_{s=1}^{t} d(u_s, v_s) = d(u^{t-1}, v^{t-1}) + d(u_t, v_t)$ and $\Pi$ is the feasible set.

Let $\pi^*(u^{t-1}, v^{t-1})$ be the transportation plan optimal for the nested distance $dl_{t-1}(\mathbb{P}, \widetilde{\mathbb{P}})$

and $\pi^*(u_t, v_t|u^{t-1}, v^{t-1})$ be the transportation plan optimal for the Kantorovich-Wasserstein distance $d_{KA}\big(P_t(\cdot|u^{t-1}), \widetilde{P}_t(\cdot|v^{t-1})\big)$. Using the transportation plan $\pi(du^t, dv^t) = \pi^*(du_t, dv_t|u^{t-1}, v^{t-1})\pi^*(du^{t-1}, dv^{t-1}) \in \Pi$ we bound the nested distance as

$$\begin{aligned} dl_t(\mathbb{P}, \widetilde{\mathbb{P}}) \leq & \int [d(u^{t-1}, v^{t-1}) + d(u_t, v_t)]\pi^*(du_t, dv_t|u^{t-1}, v^{t-1})\pi^*(du^{t-1}, dv^{t-1}) = \\ & = dl_{t-1}(\mathbb{P}, \widetilde{\mathbb{P}}) + \int d_{KA}\big(P_t(\cdot|u^{t-1}), \widetilde{P}_t(\cdot|v^{t-1})\big)\pi^*(du^{t-1}, dv^{t-1}), \end{aligned}$$

where $d_{KA}\big(P_t(\cdot|u^{t-1}), \widetilde{P}_t(\cdot|v^{t-1})\big) \leq \sup_{u^{t-1}, v^{t-1}} d_{KA}\big(P_t(\cdot|u^{t-1}), \widetilde{P}_t(\cdot|v^{t-1})\big)$ and, hence,

$$dl_t(\mathbb{P}, \widetilde{\mathbb{P}}) \leq dl_{t-1}(\mathbb{P}, \widetilde{\mathbb{P}}) + \sup_{u^{t-1}, v^{t-1}} d_{KA}\big(P_t(\cdot|u^{t-1}), \widetilde{P}_t(\cdot|v^{t-1})\big),\ \forall t = 1, ..., T.$$

A similar inequality can be written for clairvoyant tree structures, i.e.,

$$dl_t(\mathbb{P}_{c(t)}, \widetilde{\mathbb{P}}_{c(t)}) \leq dl_{t-1}(\mathbb{P}_{c(t)}, \widetilde{\mathbb{P}}_{c(t)}) + \sup_{u^{T-2}, v^{T-2}} d_{KA}\big(P^{t:T}(\cdot|\mathscr{F}_T(u^{T-2})), \widetilde{P}^{t:T}(\cdot|\mathscr{F}_T(v^{T-2}))\big). \quad (13)$$

Importantly, the distributions $P^{t:T}(\cdot|\mathscr{F}_T(u^{T-2}))$ and $\widetilde{P}^{t:T}(\cdot|\mathscr{F}_T(v^{T-2})$ are conditional on the filtration structure $\mathscr{F}_T$, which is different from the case with non-clairvoyant nested distributions and the bound (6). Further, note that $dl_{t-1}(\mathbb{P}, \widetilde{\mathbb{P}}) = dl_{t-1}(\mathbb{P}_{c(t)}, \widetilde{\mathbb{P}}_{c(t)})$ due to the fact that the final stage of these trees is $t-1$. From this follows the set of inequalities

$$\begin{aligned} & dl_T(\mathbb{P}_{c(T)}, \widetilde{\mathbb{P}}_{c(T)}) \leq dl_{T-1}(\mathbb{P}, \widetilde{\mathbb{P}}) + \sup_{u^{T-1}, v^{T-1}} d_{KA}\big(P_T(\cdot|\mathscr{F}_T(u^{T-1})), \widetilde{P}_T(\cdot|\mathscr{F}_T(v^{T-1}))\big), \\ & \vdots \\ & dl_t(\mathbb{P}_{c(t)}, \widetilde{\mathbb{P}}_{c(t)}) \leq dl_{t-1}(\mathbb{P}, \widetilde{\mathbb{P}}) + \sup_{u^{T-2}, v^{T-2}} d_{KA}\big(P^{t:T}(\cdot|\mathscr{F}_T(u^{T-2})), \widetilde{P}^{t:T}(\cdot|\mathscr{F}_T(v^{T-2}))\big), \\ & dl_{t-1}(\mathbb{P}_{c(t-1)}, \widetilde{\mathbb{P}}_{c(t-1)}) \leq dl_{t-2}(\mathbb{P}, \widetilde{\mathbb{P}}) + \sup_{u^{T-3}, v^{T-3}} d_{KA}\big(P^{t-1:T}(\cdot|\mathscr{F}_T(u^{T-3})), \widetilde{P}^{t-1:T}(\cdot|\mathscr{F}_T(v^{T-3}))\big). \end{aligned}$$

Based on the obtained recursive structure and using the fact that $dl_T(\mathbb{P}_{c(T)}, \widetilde{\mathbb{P}}_{c(T)}) = dl_T(\mathbb{P}, \widetilde{\mathbb{P}})$, we claim the statement of the corollary. Further, by the triangle inequality and given the Lipschitz property with a constant $K_{c(t)}$ for the distribution $P^{t:T}(\cdot|\mathscr{F}_T(u^{t-1}))$, we claim

$$\begin{aligned} & d_{KA}\big(P^{t:T}(\cdot|\mathscr{F}_T(u^{t-1})), \widetilde{P}^{t:T}(\cdot|\mathscr{F}_T(v^{t-1}))\big) \leq \\ & \leq d_{KA}\big(P^{t:T}(\cdot|\mathscr{F}_T(u^{t-1})), P^{t:T}(\cdot|\mathscr{F}_T(v^{t-1}))\big) + d_{KA}\big(P^{t:T}(\cdot|\mathscr{F}_T(v^{t-1})), \widetilde{P}^{t:T}(\cdot|\mathscr{F}_T(v^{t-1}))\big) \leq \\ & \leq K_{c(t)} d(u^{t-1}, v^{t-1}) + \sup_{v^{t-1}} d_{KA}\big(P^{t:T}(\cdot|\mathscr{F}_T(v^{t-1})), \widetilde{P}^{t:T}(\cdot|\mathscr{F}_T(v^{t-1}))\big) = \\ & = K_{c(t)} d(u^{t-1}, v^{t-1}) + d_{KA}(P^{t:T}, \widetilde{P}^{t:T}) \end{aligned}$$

and, thus, $dl_t(\mathbb{P}_{c(t)}, \widetilde{\mathbb{P}}_{c(t)}) \leq (1 + K_{c(t)}) dl_{t-1}(\mathbb{P}, \widetilde{\mathbb{P}}) + d_{KA}(P^{t:T}, \widetilde{P}^{t:T})$, $\forall t = 1, ..., T$. Applying this recursion for all $t$ we obtain the bound similar to the bound (7), but with joint distributions and filtration $\mathscr{F}_T$ at stages $t, ..., T$:

$$\inf_{\widetilde{\mathbb{P}}} dl(\mathbb{P}, \widetilde{\mathbb{P}}) \leq \sum_{t=1}^{T} \inf_{\widetilde{P}^{t:T}} d_{KA}(P^{t:T}, \widetilde{P}^{t:T}) \prod_{s=t+1}^{T} (K_{c(t)} + 1), \quad (14)$$

where $d_{KA}(P^{t:T}, \widetilde{P}^{t:T}) = \sup_{v^{t-1}} d_{KA}\big(P^{t:T}(\cdot|\mathscr{F}_T(v^{t-1})), \widetilde{P}^{t:T}(\cdot|\mathscr{F}_T(v^{t-1}))\big)$. □

Note that it is sufficient to use the upper bound (6) if the components $\xi_t$, $\forall t = 1, ..., T$ of the stochastic process $\xi = (\xi_1, ..., \xi_T)$ are independent and identically distrubuted (i.i.d.) random variables: this is due to the implied independence between time stages. However, if the components $\xi_t$, $\forall t = 1, ..., T$ of the stochastic process $\xi = (\xi_1, ..., \xi_T)$ are interdependent on one another, the upper bounds (9) and (12) become crucial.

## 3. Multi-Stage Scenario Quantization

The problem to find a scenario tree, which is closest to a given continuous-state stochastic process in the sense of the minimal nested distance, given the total number of nodes only, seems to be too complex to be solved in a satisfactory manner. However, if the topology of the tree is chosen beforehand and fixed, the minimal distance can be bounded by a computable expression. To explain this, suppose that the scenario tree topology is fixed with a total of $1 + N_1 + N_2 + \cdots + N_{T-1} + N_T$ nodes, where $N_T = N$[1]:

$$\begin{array}{ll} (0,1) & \text{at stage } 0 \ (\textit{the root}) \\ (1,1),\ldots,(1,N_1) & \text{at stage } 1 \\ (2,1),\ldots,(2,N_2) & \text{at stage } 2 \\ \vdots & \\ (T-1,1),\ldots,(T-1,N_{T-1}) & \text{at stage } T-1 \\ 1,\ldots,N & \text{at stage } T \ (\textit{the leaves}). \end{array}$$

For each stage $t$, there are predecessor functions $pred_t(i)$, $t = 1, \ldots, T-1$, $i = 1, \ldots, N$ such that $(t, pred_t(i))$ is the predeccessor of scenario $i$ (i.e., leaf $i$) at stage $t$. Clearly, $pred_0(i) = 1 \ \forall i$.

The value associated with the node $(t,i)$ is denoted by $z_t^{(i)}$. Node values can be scalars or vectors. Let $\mathscr{Z}$ be the vector space of all $[N \times T]$ matrices of the form

$$Z = \begin{pmatrix} z_1^{(1)} & \cdots & z_{T-1}^{(1)} & z_T^{(1)} \\ \vdots & \ddots & \vdots & \vdots \\ z_1^{(N)} & \cdots & z_{T-1}^{(N)} & z_T^{(N)} \end{pmatrix}$$

and let $\mathscr{Z}^{na}$ be the linear subspace of $\mathscr{Z}$, where the following equations hold

$$\mathscr{Z}^{na} = \left\{ Z \in \mathscr{Z} : z_t^{(i)} = z_t^{(j)} \text{ if } pred_t(i) = pred_t(j) \right\}.$$

We refer to the linear subspace $\mathscr{Z}^{na}$ of $\mathscr{Z}$ as the space of *(nonanticipative) node values*. If $Z \in \mathscr{Z}^{na}$, the $z$-values can be assigned in a canonical manner to the nodes of the tree with the given topology.

Suppose now that the tree topology is fixed and the values at the nodes are given by $Z \in \mathscr{Z}^{na}$. The node probabilities may vary. Denoting the set of all the nested tree distributions by $\tilde{\mathscr{P}}$, it turns out that it is easier to find the minimum distance w.r.t the set $\tilde{\mathscr{P}}$ than to find the distance w.r.t. a fixed $\tilde{\mathbb{P}} \in \tilde{\mathscr{P}}$.

To express this minimal distance, consider a stochastic process $(\xi_1, \ldots, \xi_T)$ together with the filtration generated by this process. Denote by $\mathbb{P}$ the nested distribution of this process. Then the minmal nested distance between $\mathbb{P}$ and the set $\tilde{\mathscr{P}}$ with fixed node values $Z \in \mathscr{Z}^{na}$ is given by the value of the following expectation:

$$\mathbb{E}_\xi \left[ \min_i \sum_{t=1}^{T} \|z_t^{(i)} - \xi_t\| : i = 1, \ldots, N \right].$$

[1] Note that the leaf nodes are denoted simply by $1, \ldots, N$, where $N$ is the number of scenarios.

The more involved problem of finding the nearest (in the sense of the minimal nested distance) tree with a given topology and node values $Z$ as decision variables is the global solution of the following stochastic optimization problem:

$$\min_{Z \in \mathscr{Z}^{na}} \mathbb{E}_\xi \left[ \min_i \sum_{t=1}^{T} \|z_t^{(i)} - \xi_t\| : i = 1, \ldots, N \right]. \tag{15}$$

Unfortunately, the problem (15) is nonconvex and any computed local minimizer is just an upper bound of the minimal nested distance. Nevertheless, any algorithm which reduces the objective value in (15) gives an improvement in nested distance.

Evidently, a *stochastic projected gradient algorithm* would be applicable for the problem (15). However, due to the high nonconvexity with very many local minima it is not advisable to start with such a procedure. In this article, we propose to construct the scenario tree via a two-step procedure combining *a forward and a backward step*: in the forward step, we use the discrete approximations of the stage-wise conditional distributions, which guarantee an upper bound for the nested distance. In the backward step, we improve the distance by applying a stochastic projected gradient algorithm to solve a problem of type (15), but only for two subsequent stages in line with upper bounds (10) and (11). Before we proceed with the description of the algorithm, we discuss methods, which can be used in the forward step.

### 3.1. Existing scenario quantization methods

Current literature relies on a stage-wise discretization of continuous distributions $P_t$, $\forall t = 1, ..., T$ and provides an approximation to the upper bound (6), but not its minimum [36, 49, 59]. In particular, multiple methods exist for the approximation of a random variable $\xi_t$ by a measurable function sitting on $N$ supporting points, e.g., statistical methods described in the works [9, 27, 30], as well as the following general frameworks:

***Monte-Carlo (random) generation*** randomly selects $N$ points from the distribution function $P_t$ and assigns equal probabilities to them [16, 59, 49];

***Quasi Monte-Carlo generation*** replaces random samples of the Monte-Carlo method by deterministic points that are uniformly distributed in $[0,1]^D$ [59, 49];

***Stage-wise optimal quantization*** minimizes the Kantorovich-Wasserstein distance $d_{KA}(P_t, \sum_{i=1}^{N} p_i \delta_{z_t^{(i)}})$ between the continuous distribution function $P_t$ and its discrete approximation $\sum_{i=1}^{N} p_i \delta_{z_t^{(i)}}$. The optimal supporting points $z_t^{(i)}$, $i = 1, ..., N_t$ are computed by the minimization of a functional $\int \min_i d(x, z_t^{(i)}) P_t(dx)$ over the decision variables $z_t^{(1)}, ..., z_t^{(N_t)}$ [17, 33, 36, 53, 55, 68, 69]. Discretizing stage-wise distributions $P_t$, $\forall t = 1, ..., T$ given the number $N_t$ of supporting points corresponding to the structure of the scenario tree (see [63, 65]), we obtain a discrete approximation of the nested distribution $\mathbb{P}$, its corresponding multivariate probability $P^{1:T}$ and a chain of continuous conditional distributions (see [63, 65]).

Figure 2 compares the Monte-Carlo sampling (Figure 2a) with the stage-wise optimal quantization (Figure 2b) for a one-dimensional Gaussian distribution sitting on the tree with a given structure (i.e., $N = 5$ optimal supporting points for each conditional distribution). Here, we discretize each of the conditional distributions at stages $t = 1, 2, 3$ with 5 supporting points.

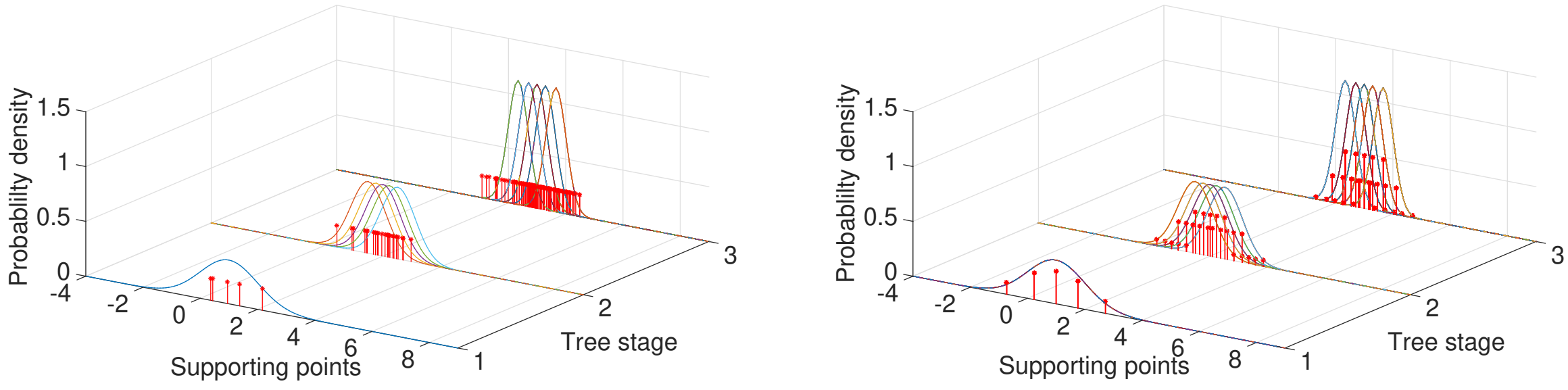

a. Random generation. b. Stage-wise optimal quantization.

**Figure 2 Optimal quantization and random generation for a Gaussian distribution with mean** $\mu = (1,\ 2,\ 3)$ **and covariance** $C = (1,\ 0.3,\ 0;\ 0.3,\ 0.7,\ 0.3;\ 0,\ 0.3,\ 0.5)$**.**

As every quantization method needs to reduce the nested distance (5) to guarantee a fine approximation quality in multi-stage stochastic optimization, the *size* and *speed* of such a reduction become the most important factors in evaluating the performance of the method. To observe the reduction of the nested distance between the nested distribution $\mathbb{P}$ and its approximation $\widetilde{\mathbb{P}}^+$ given the methods in Figure 2, we consider a sequence of scenario trees $\mathbb{Q}_b$ with the same number of supporting points $b$ for each conditional distribution. The number $b$ is defined as the *branchiness factor* of a scenario tree in the work [63]. The convergence of the nested distance is guaranteed in probability in line with the triangle inequality $dl(\mathbb{P}, \widetilde{\mathbb{P}}^+) \leq dl(\mathbb{P}, \mathbb{Q}_b) + dl(\widetilde{\mathbb{P}}^+, \mathbb{Q}_b)$, where the *error term* goes to zero (i.e., $dl(\mathbb{P}, \mathbb{Q}_b) \to 0$) and the *main term* converges to the true nested distance (i.e., $dl(\widetilde{\mathbb{P}}^+, \mathbb{Q}_b) \to dl(\mathbb{P}, \widetilde{\mathbb{P}}^+)$) as the branchiness factor increases (i.e., $b \to \infty$) independently of the quantization method of the tree $\mathbb{Q}_b$. This is due to Lemmas 2-3 for a sequence of scenario trees $\mathbb{Q}_b$ with increasing branchiness (see Appendix 6). Figures 3 and 4 demonstrate that the nested distance is lower if $\widetilde{\mathbb{P}}^+$ is stage-wise optimal. Thus, the stage-wise optimal quantization provides a finer approximation of the nested distribution $\mathbb{P}$ than the random (Monte-Carlo) sampling.

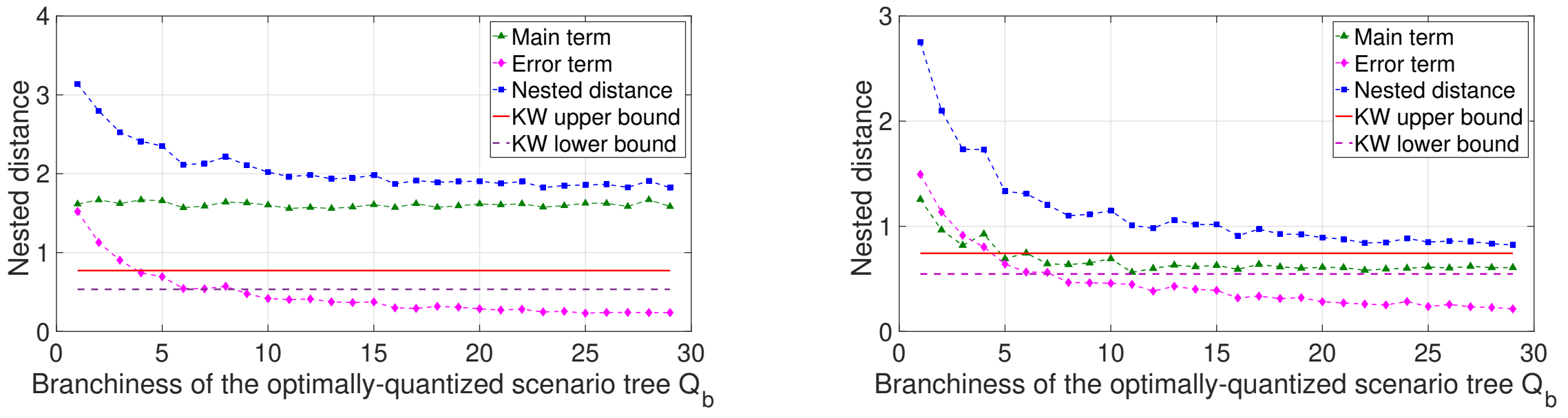

a. Randomly sampled $\widetilde{\mathbb{P}}^+$. b. Optimally-quantized $\widetilde{\mathbb{P}}^+$.

**Figure 3 Nested distance between Gaussian stochastic process and its discrete approximation.**

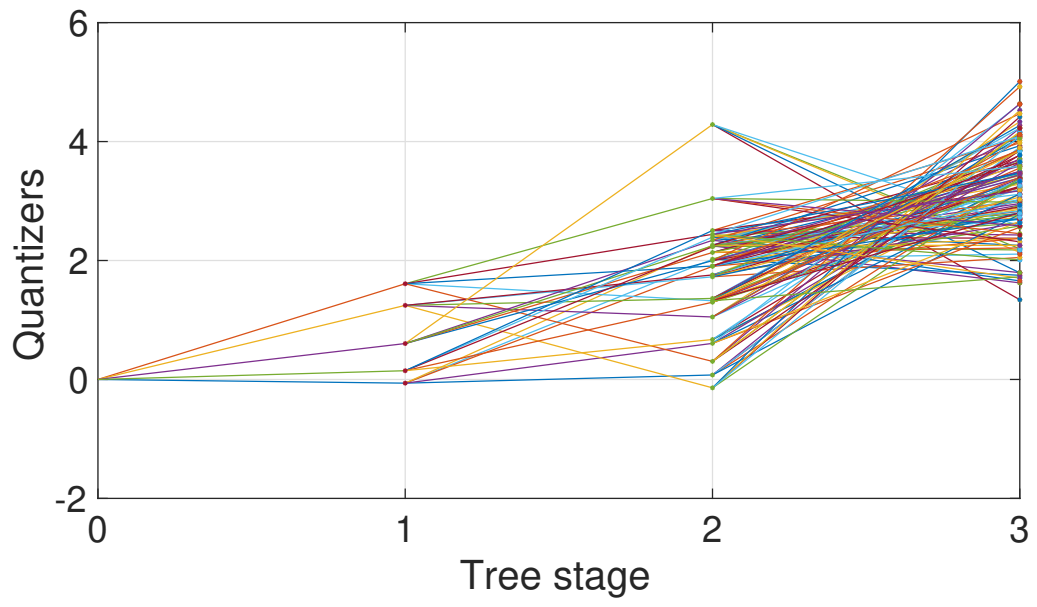

a. Randomly sampled $\widetilde{\mathbb{P}}^{+}$.

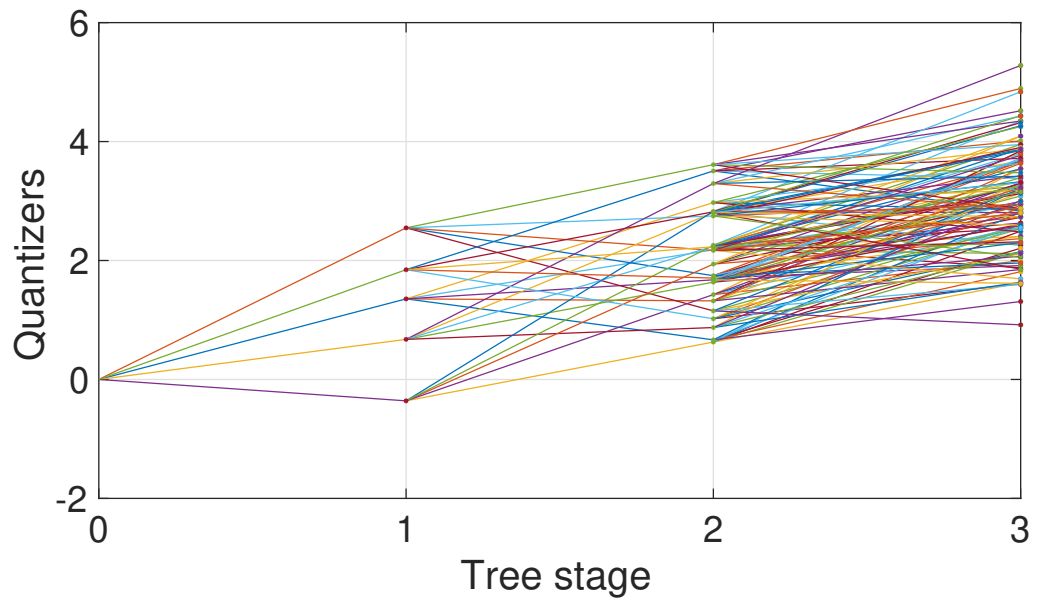

b. Optimally-quantized $\widetilde{\mathbb{P}}^{+}$.

**Figure 4 Nested distance between Gaussian stochastic process and its discrete approximation.**

Besides the size of the nested distance reduction, quantization methods clearly differ in the number of supporting points, which they require for the same approximation accuracy in multi-stage stochastic optimization. Indeed, Monte-Carlo sampling generates points randomly and, thus, heavily relies on the Law of Large Numbers for the fine approximation. The number of supporting points required by the stage-wise optimal quantization is lower due to the optimality of the point sampling method at each tree stage.

In this article, we propose a method which guarantees a fine approximation accuracy with even less number of points. To achieve this, we consider interactions between different tree stages and aim to reduce the nested distance upper bounds (10) and (11), tightly bounding the optimality gap between the initial and the approximate problems.

Further, we refer to the discrete nested distribution $\widetilde{\mathbb{P}}^{+}$ and its chain of conditional probabilities $\widetilde{P}_1^{+}$, $\widetilde{P}_2^{+}$,...,$\widetilde{P}_T^{+}$ as distributions quantized using the stage-wise optimal quantization. Though the stage-wise optimal quantization provides a finer approximation in comparison to the random sampling (Figure 3), it does not necessarily lead to the minimization of the upper bound (6) as the realization paths $u^{t-1}$ and $v^{t-1}$ are not necessarily optimal (i.e., they do not necessarily provide the supremum of the Kantorovich-Wasserstein distance at the stage $t$ of the formula (13)). Clearly, the stage-wise optimal quantization also cannot guarantee optimality for the nested distance between the stochastic process and a tree, i.e.,

$$\widetilde{\mathbb{P}}^{+} \neq \mathbb{P}^{*} \in \arg\min_{\widetilde{\mathbb{P}}} dl(\mathbb{P}, \widetilde{\mathbb{P}}),$$

where we denote by $\mathbb{P}^{*}$ the optimal nested distribution in the sense of minimal nested distance.

### 3.2. Forward-Backward Quantization

In order to improve the nested distribution approximation $\widetilde{\mathbb{P}}^{+}$, we account for interstage dependencies using joint distributions in line with the upper bound (12). For this, we note that the $T$-th term of the upper bound (12) coincides with the corresponding term of the bound (6) and, thus, distribution quantizations at the stage $T$ are equivalent for these bounds. Nevertheless, the quantization implied by the upper bound (12) differs from the one implied by the bound (6) for stages $t = 1,...,T-1$ as the bound (12) is conditional on the feasible future scenarios.

Overall, we develop a novel scenario quantization method on the basis of the upper bound (12) and numerically demonstrate that its minimization results in a finer approximation of the process $\mathbb{P}$. For this, we propose a two-step procedure, which initially computes the stage-wise optimal quantizers and, afterwards, improves them running the procedure backwards from the leaves of the scenario tree (i.e., starting at stage T). We obtain the improved quantizers at stage $t$ using the filtration structure of stage $T$.

Before we proceed with the description of the algorithm, we provide an explanatory example (see Example 1 and Appendix 6.1).

EXAMPLE 1. *Consider a two-stage scenario tree in Figure 5a, $t = 1,2$. This tree approximates a continuous-state stochastic process $\xi = (\xi_1, \xi_2)$, where one-dimensional random variable $\xi_2$ depends on the realization of the one-dimensional random variable $\xi_1$. The initial discrete quantization of the tree is shown in Figure 5b, where the quantizers are denoted by $(z^{(1)}, z^{(2)}, ..., z^{(10)})'$.*

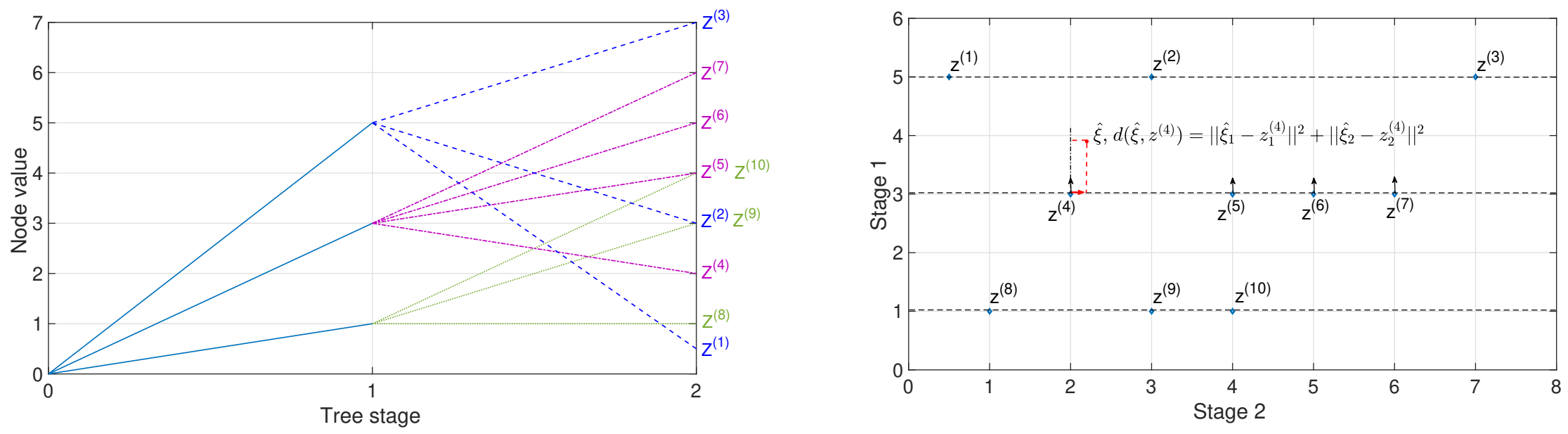

a. Scenario tree. b. Tree quantizers.

**Figure 5 Nested distance minimization given the filtration structure of a two-stage scenario tree.**

*Considering the $l_2$-distance metric $d(\cdot,\cdot)$ and $r = 2$ for the nested distance minimization (see Definition 1), we obtain the following optimization problem:*

$$\begin{aligned} \underset{z_i^{(1)}, z_i^{(2)}}{\textit{min}} \quad & \mathbb{E}\left[\min_i \left((z_1^{(i)} - \xi_1)^2 + (z_2^{(i)} - \xi_2)^2\right)\right] \\ \textit{subject to} \quad & z_1^{(1)} = z_1^{(2)} = z_1^{(3)}, \\ & z_1^{(4)} = z_1^{(5)} = z_1^{(6)} = z_1^{(7)}, \\ & z_1^{(8)} = z_1^{(9)} = z_1^{(10)}. \end{aligned} \tag{16}$$

*The constraints of the problem (16) are convex and can be represented as a set $\Xi = \{z \in \mathbb{R}^{10} : Az = 0\}$ with a given matrix A. These constraints set first coordinates of vectors $z^{(i)}$ equal to each other, following the structure of the scenario tree in Figure 5a. We denote by $n_i$ the corresponding number of equal first coordinates for each quantizer i (i.e., $n_1 = n_2 = n_3 = 3$, $n_4 = n_5 = n_6 = n_7 = 4$, $n_8 = n_9 = n_{10} = 3$). Clearly, if the height of the scenario tree would be greater than two, one could generalize this notation to $n_t^{(i)}$ specifying the stage t, at which the quantizers are equal. We follow this notation in the proposed Algorithm 1. Further note that $z_1^{(i)}$ and $z_2^{(i)}$ represent multi-dimensional vectors in the general case (Algorithm 1).*

*The optimization problem (16) can be solved via stochastic iterations of the* ***projected gradient method*** *with a step-size $t_k$ (see Appendix 6.1). For this, we randomly sample a point $\hat{\xi}(k) = (\hat{\xi}_1(k), \hat{\xi}_2(k))$ at each iteration $k$ following the stochastic process $\xi$ and find the quantizer $z^{(i)}$ such that $i \in$* $\arg\min_i \left( (z_1^{(i)} - \hat{\xi}_1(k))^2 + (z_2^{(i)} - \hat{\xi}_2(k))^2 \right)$*. Afterwards, we update the coordinates as follows:*

$$z(k+1) = P_{\Xi}\left(z(k) - t_k \nabla(z(k), \hat{\xi}(k))\right) = z(k) - t_k P_{\Xi} \nabla(z(k), \hat{\xi}(k)),$$

*where $P_{\Xi} = I - A'(AA')^{-1}A$ is the projection onto the convex set of constraints $\Xi$ and $\nabla(z(k), \hat{\xi}(k))$ is the* stochastic gradient function *(see Appendix 6.1). Overall, the iterations can be written in the following form:*

$$\begin{cases} z_1^{(i)}(k+1) = z_1^{(i)}(k) - \dfrac{2t_k}{n_i}(z_1^{(i)}(k) - \hat{\xi}_1(k)), \\ z_1^{(j)}(k+1) = z_1^{(j)}(k) - \dfrac{2t_k}{n_i}(z_1^{(j)}(k) - \hat{\xi}_1(k)), \; \forall j : z_1^{(i)}(k) = z_1^{(j)}(k), \\ z_2^{(i)}(k+1) = z_2^{(i)}(k) - 2t_k(z_2^{(i)}(k) - \hat{\xi}_2(k)), \end{cases}$$

*where the division by $n_i$ arises because of the projection of the gradient vector onto the convex set of constraints of the problem (16). Clearly, the gradient differs dependent on the dimensionality of the stochastic process, the chosen distance metric $d(\cdot,\cdot)$ and the chosen parameter $r$ in the nested distance minimization problem in Definition 1. For example, given a multi-dimensional process $\xi = (\xi_1, \xi_2)$, $r = 1$ and the $l_2$-distance metric, one would obtain the following iterations:*

$$\begin{cases} z_1^{(i)}(k+1) = z_1^{(i)}(k) - \dfrac{t_k}{n_i} \dfrac{z_1^{(i)}(k) - \hat{\xi}_1(k)}{\|z_1^{(i)}(k) - \hat{\xi}_1(k)\|_2}, \\ z_1^{(j)}(k+1) = z_1^{(j)}(k) - \dfrac{t_k}{n_i} \dfrac{z_1^{(j)}(k) - \hat{\xi}_1(k)}{\|z_1^{(j)}(k) - \hat{\xi}_1(k)\|_2}, \; \forall j : z_1^{(i)}(k) = z_1^{(j)}(k), \\ z_2^{(i)}(k+1) = z_2^{(i)}(k) - t_k \dfrac{z_2^{(i)}(k) - \hat{\xi}_2(k)}{\|z_2^{(i)}(k) - \hat{\xi}_2(k)\|_2}. \end{cases}$$

*Here, the sequence of step sizes $t_k$ can be chosen in line with the* ***square summable but not summable*** *rule satisfying $\sum_{k\to\infty} t_k^2 < \infty$, $\sum_{k=1}^{\infty} t_k = \infty$. Further, one can rely on the* ***Polyak's step size*** *of the type $t_k = \frac{1}{k\|\nabla(z(k),\hat{\xi}(k))\|_2^2}$. In the proposed Algorithm 1 we use the* ***diminishing step size*** *such that $\sum_k t_k = \infty$ and $\sum_n t_k^2 < \infty$. The convergence of the stochastic gradient algorithm was studied by numerous authors, see e.g., [42, 54, 58].* □

Note that Example 1 corresponds to the local nested distance minimizer, where the optimal plan $\pi$ transports the whole mass corresponding to the point $\xi$. This is because the stochastic process $\xi$ is continuous-state, which implies point-to-point transportation (see Definition 1). Oppositely, if both the initial and the approximate processes were discrete, the optimal transportation plan would need to be computed as well. In this article, we focus on the approximation of continuous-state stochastic processes by discrete scenario trees and, thus, on point-to-point transportations.

Below, we describe Algorithm 1 for the case with the $l_2$-distance measure $d(\cdot,\cdot)$ and $r = 2$. Overall, the method can be subdivided into a **forward** and a **backward** step.

***Forward step:*** Stage-wise optimal quantizers approximate the stochastic process. We denote the scenario process between stages $t$ and $T$ by $z^{t:T,(i)} = (z_t^{(i)}, ..., z_T^{(i)})$, $\forall t, i$, where $t$ is the time stage;

***Backward step:*** We adapt quantizers $z_t^{(i)}$, $\forall t, i$ so that the Kantorovich-Wasserstein distance between the joint probability distribution of stages $t, ..., T$ and its discrete approximation is minimized under a fixing constraint for the scenario tree filtration structure (see Algorithm 1).

**Algorithm 1** Forward-Backward Optimal Quantization ($l_2$-distance measure, $r = 2$).

Compute stage-wise optimal quantizers $z_t^{(i)}$, $\forall t, i$.

**for** $t = T-1, ..., 1$ **do**

Let $I_t$ be the number of nodes at the stage $t-1$ of the scenario tree. By this, $I_t$ represents the number of subtrees with roots at stage $t-1$ and, equivalently, the number of conditional distributions at stage $t$. We assign indices $I = 1, ..., I_t$ to these subtrees and note that each subtree $I$ has multiple scenarios associated with it. Corresponding sets of scenario indices are denoted by $S_I$, $\forall I$.

**for** $I = 1, ..., I_t$ **do**

Sample $K$ random points $\hat{\xi}(k) = (\hat{\xi}_t(k), \hat{\xi}_{t+1}(k), ..., \hat{\xi}_T(k))'$, $k = 1, ..., K$ from the joint distribution at the $I$-th subtree. For each point $k$, find the nearest scenario $z^{t:T,(i)}(k) = (z_t^{(i)}(k),\ z_{t+1}^{(i)}(k), ..., z_T^{(i)}(k))'$, i.e., $i \in \arg\min\{d(\hat{\xi}(k), z^{(i)}(k)),\ i \in S_I\}$, where $z^{t:T,(i)}(1) = z^{t:T,(i)}$ at the first iteration. As different scenarios may be equal at some time stages, we let $n_t^{(i)}$ be the number of scenarios with equal $t$-th coordinate $z_t^{(i)}(k)$ at any stage $t$ of a particular subtree (see Figure 6b).

We update the quantizers as follows:

$$\begin{cases} z_t^{(j)}(k+1) = z_t^{(j)}(k) - \dfrac{1}{n_t^{(i)} k^{0.75}} \left( z_t^{(j)}(k) - \hat{\xi}_t(k) \right), \ \forall j \in S_I : z_t^{(j)}(k) = z_t^{(i)}(k). \\ z_{t+1}^{(m)}(k+1) = z_{t+1}^{(m)}(k) - \dfrac{1}{n_{t+1}^{(i)} k^{0.75}} \left( z_{t+1}^{(m)}(k) - \hat{\xi}_{t+1}(k) \right), \ \forall m \in S_I : z_{t+1}^{(m)}(k) = z_{t+1}^{(i)}(k), \\ \dots \\ z_T^{(i)}(k+1) = z_T^{(i)}(k) - \dfrac{1}{n_T^{(i)} k^{0.75}} \left( z_T^{(i)}(k) - \hat{\xi}_T(k) \right). \end{cases}$$

To estimate the probabilities of scenarios, we compute the number $K_t^{(i)}$ of random samples $\hat{\xi}(k)$ nearest to the point $(z_t^{(i)}(k),\ z_{t+1}^{(i)}(k), ..., z_T^{(i)}(k))'$. The probabilities are computed as $p_t^{(i)} = \frac{K_t^{(i)}}{K}$, $\forall t, i$.

**end for**

**end for**

In Algorithm 1, we use the diminishing step size $t_k = \frac{n_t^{(i)}}{2k^{0.75}}$. However, the use of other step sizes (e.g., square summable but not summable or Polyak's step sizes) is possible in stochastic gradient methods [42, 54, 58].

Using the well-known facility location problem as an analogy, note that Algorithm 1 decides where to build stores in a city accounting for the distribution of its population and having less degrees of freedom, e.g., moving the stores along given streets only. Overall, Algorithm 1 is a *learning*-type algorithm, which minimizes the nested distance upper bound via varying stage-wise quantizers in a joint setup with a fixed filtration structure. It aims to reduce the nested distance taking all scenario tree stages into account. This is unlike the stage-wise optimal quantization, which does not account for interdependencies between current and future stages. Note that the iteration of the algorithm updates the quantizer at stage $t$ and all the future quantizers (see Figure 6).

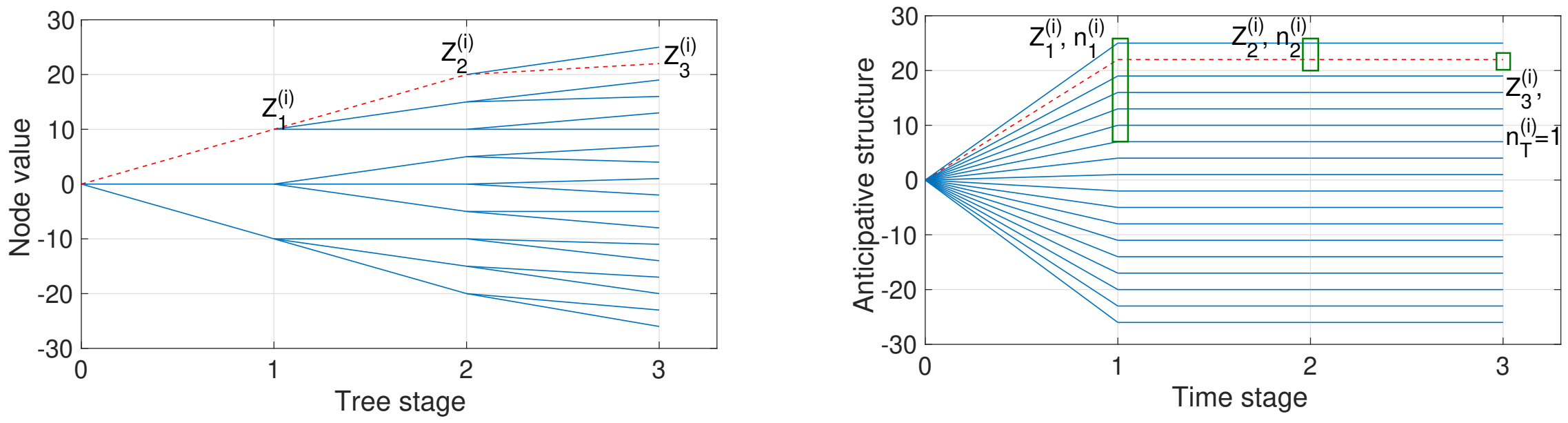

a. Initial tree. b. Anticipative tree.

**Figure 6 Filtration structure of a multi-stage scenario tree.**

In the numerical part of the article, we test Algorithm 1. Dependent on the chosen upper bound, one can avoid the adaptation of some (or even all) of the future scenarios. For example, working with the upper bounds (10) and (11), one adapts the quantizers at stages $t$ and $t+1$ $\forall t$ at each iteration $k$, which makes the algorithm more efficient. By this, one receives at least an upper bound on the two-stage nested distances, while the global solution would lead to the minimum of the nested distance itself. In Sections 3.3 and 4 we test an efficient adaptation of the Algorithm 1 avoiding the change of all the future scenarios $t+1,\ldots,T$. Figure 7a schematically demonstrates the step of this method using an example of a bivariate Gaussian distribution between random variables $\xi_t$ and $\xi_{t+1}$. In Figure 7a, we adapt the position of quantizers which belong to the stage $t$ requiring the non-anticipativity necessary for the upper bound. Oppositely, relaxing the tree structure as shown in Figure 7b results in a lower bound for the nested distance.

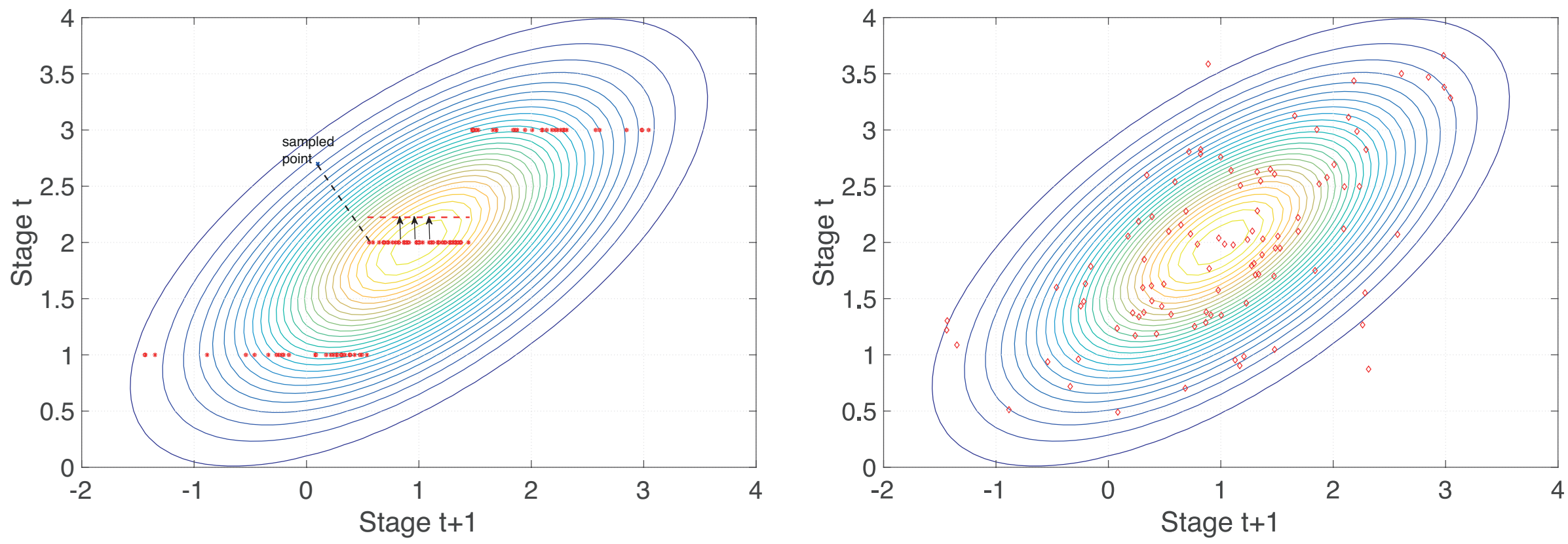

a. Random samples with a tree structure. b. Random samples without a tree structure.

**Figure 7 Two-stage quantization with and without a tree structure.**

### 3.3. Functional properties

To analyze the performance of Algorithm 1, we consider a stochastic process $(\xi_1,\dots,\xi_T)$ following a multivariate Gaussian distribution with mean vector $\mu=(\mu_1,\dots,\mu_T)$ and a non-singular covariance matrix $C=(c_{s,t})_{s,t=1,\dots,T}$. The conditional expectation and the variance of the random variable $\xi_t$ given its historical process realization $\xi^{t-1}$ are denoted by $\mu_t(\xi^{t-1})$ and $\sigma_t^2(\xi^{t-1})$, respectively, in line with the work [34]. The conditional distribution of $\xi_t$ given $\xi^{t-1}$ can be written in the closed form, i.e.,

$$\xi_t \mid \xi^{t-1} \sim \mathcal{N}\left(\mu_t+(\xi^{t-1}-\mu^{t-1})C_{t-1}^{-1}c^t, c_{t,t}-(c^t)^\top C_{t-1}^{-1}c^t\right), \tag{17}$$

where $\mu^{t-1}=\left(\mu_1,...,\mu_{t-1}\right)$ is the mean process up to time $t-1$ and $C_t$ is the covariance matrix dissected into $C_t=\left(\begin{array}{c|c} C_{t-1} & c^t \\ \hline (c^t)^\top & c_{t,t}\end{array}\right)$. For the distribution (17), the Lipschitz property holds with constants $K_t=||C_{t-1}^{-1}c^t||$, $\forall t$.

Firstly, using the Gaussian distribution, we demonstrate that an efficient version of Algorithm 1, which avoids the adaptation of quantizers at stages $t+1,...,T$, performs better than the same algorithm adapting all the quantizers at stage $t,...T$. This implies a lower upper bound value of Algorithm 1 with less quantizer updates and a non-negative difference in Figure 8. For the computation, we use the upper bound (7) evaluated at the corresponding quantizers and plot the difference in the upper bound estimates in Figure 8. For this, we vary such parameters as the height $T$ of the scenario tree, its branchiness $b$, the dimensionality $D$ of the stochastic process and the dependency multiplier $\lambda$ (i.e., the factor modifying the covariance matrix to $\lambda C$). We observe that the difference in upper bounds is larger for high dependency parameters, tree height and the process dimensionality. Thus, we use the more efficient version of Algorithm 1 with less updates of coordinates in our subsequent computations.

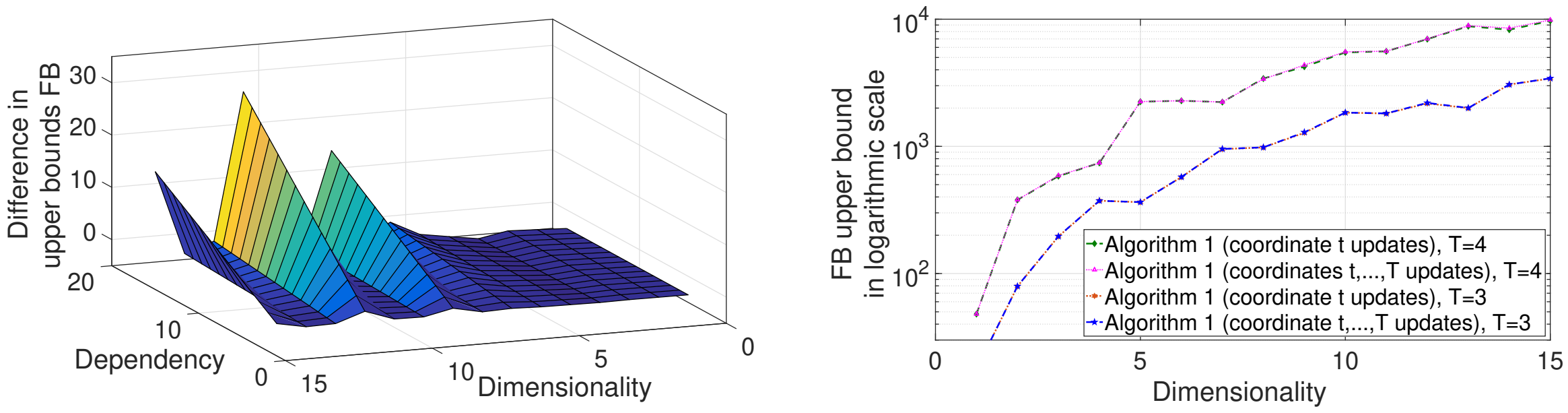

a. Tree height $T=10$. b. Dependency multiplier $\lambda=10$.

**Figure 8 Difference in upper bounds for the nested distance between stochastic process and a scenario tree.**

Figures 9-12 demonstrate that Algorithm 1 gives a better guarantee for the minimal nested distance than the stage-wise optimal quantization. The comparison is performed using the upper bound (7) evaluated at the quantizers obtained via Algorithm 1 (i.e., **FB upper bound**) and, differently, at the stage-wise optimal points (i.e., **KW upper bound**). In Figure 9a,b, we consider a binary scenario tree representing a continuous-state stochastic process $\xi$ which follows a Gaussian distribution. We observe that an increase in the process dimensionality $D$, the dependency multiplier $\lambda$, as well as in the tree height $T$ enlarges the difference between corresponding bounds.

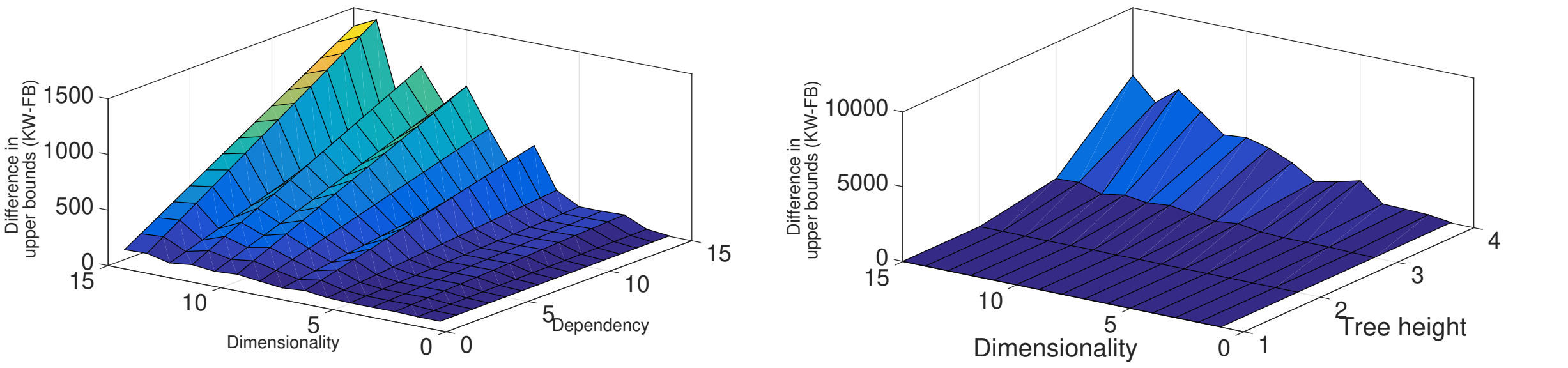

a. Tree height $T = 3$.

b. Dependency multiplier $\lambda = 10$.

**Figure 9 Difference in upper bounds for the nested distance between stochastic process and a scenario tree.**

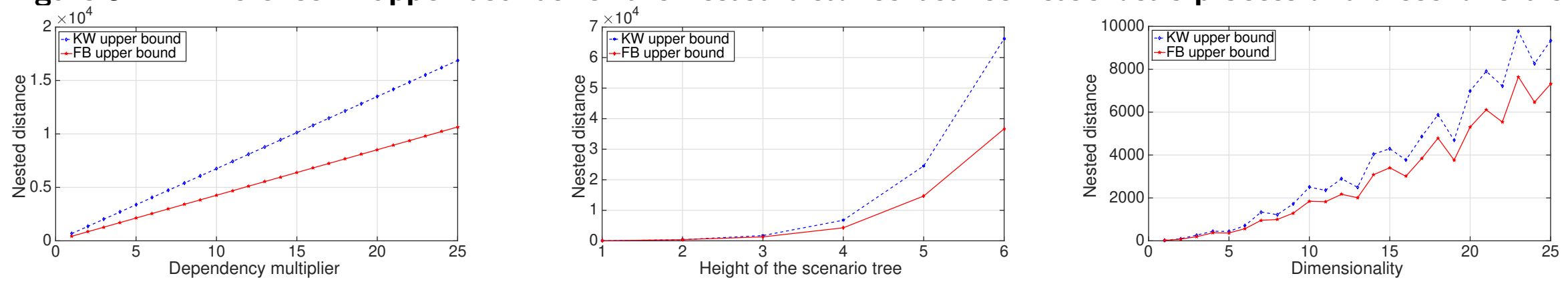

a. $T = 4, D = 9$.

c. $D = 9, \lambda = 10$

b. $T = 3, \lambda = 10$.

**Figure 10 Trends and gaps in upper bounds for the nested distance between stochastic process and a tree.**

In particular, Figure 10a shows linear increasing patterns for both upper bounds for the growing dependency multiplier, as well as the enlarging gap between the two. Similarly, the enlarging gap is observed in the exponentially increasing trends as the height $T$ of the scenario tree grows and for the fluctuating trend dependent on the dimensionality of the approximated stochastic process (Figure 10b,c).

In relation to Figure 10b, we demonstrate stage-wise terms of the upper bound (7) in Figure 11. As one can observe, some of the terms may be lower in case the stage-wise optimal quantization is used (see stage 3 in Figure 11a). Nevertheless, the sum of terms evaluated at the quantizers of Algorithm 1 stays lower than in the case with stage-wise optimal quantization: this goes in line with the idea of the nested distance minimization, which takes the dependencies between stages explicitly into account. Furthermore, Figure 11 demonstrates the behavior of the bound (7), in which joint distributions are evaluated at each stage. As expected, it imposes the weakest bound on the nested distance due to the higher optimal transportation cost.

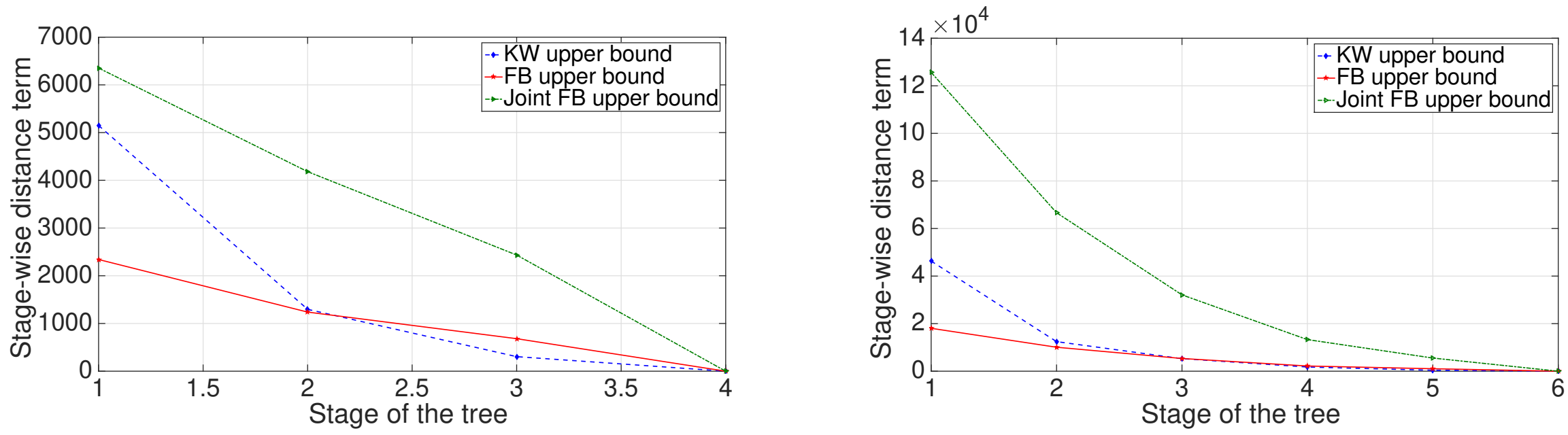

a. $D = 9, T = 4, \lambda = 10$.

b. $D = 9, T = 6, \lambda = 10$.

**Figure 11 Stage-wise terms in upper bounds for the nested distance between stochastic process and a tree.**

Next, we show the quantizers corresponding to one of the dimensions of the 9-dimensional stochastic process in Figure 12. As the FB quantization takes the dependency between present and future stages ex-plicitly into account, we observe the pronounced difference in the links of the scenario trees. At the same time, the quantizer values for this dimension are similar between two trees.

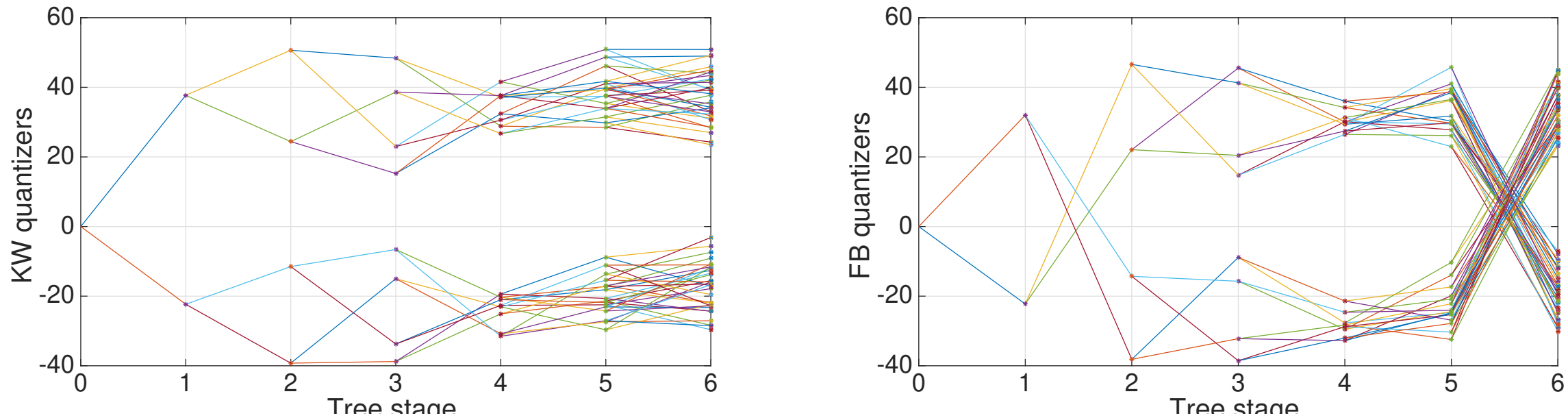

a. KW: $D = 9$, $T = 6$, $\lambda = 10$. b. FB: $D = 9$, $T = 6$, $\lambda = 10$.

**Figure 12 Difference in upper bounds for the nested distance between stochastic process and a scenario tree.**

We analyze the behavior of the nested distance estimate and compare it with its upper and lower bounds. To obtain the estimate between the continuous-state nested distribution $\mathbb{P}$ and its discrete approximation $\widetilde{\mathbb{P}}$, we use the triangle inequality $dl(\mathbb{P}, \widetilde{\mathbb{P}}) \leq dl(\mathbb{Q}_b, \widetilde{\mathbb{P}}) + dl(\mathbb{Q}_b, \mathbb{P})$ with the nested distribution $\mathbb{Q}_b$ corresponding to the tree with the branchiness factor $b$ similar to Figure 3. As shown in [63], the nested distance $dl(\mathbb{Q}_b, \widetilde{\mathbb{P}})$ converges to $dl(\mathbb{P}, \widetilde{\mathbb{P}})$ as the number of branches $b$ grows, i.e., $b \to \infty$. Therefore, Figure 13a demonstrates the convergence of the nested distance $dl(\mathbb{Q}_b, \widetilde{\mathbb{P}})$ in the number of branches $b$. Importantly, the nested distance converges to the value, which is 37.57% tighter to the lower bound than the upper bound (7) evaluated at stage-wise optimal quantizers. Additionally, Figure 13b demonstrates the convergence of the nested distance estimate for a randomly sampled tree, comparing it to just one iteration of Algorithm 1 performed on this tree. In this case, we observe ca. 13.47% outperformance of Algorithm 1 in Figure 13b.

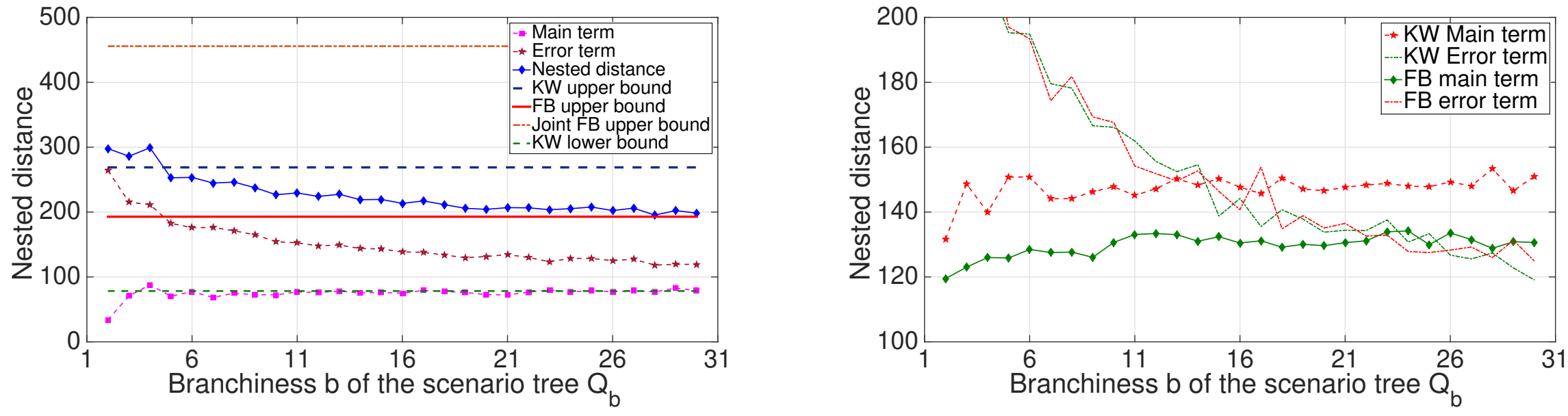

a. Nested distance bounds $D = 3$, $T = 3$, $\lambda = 10$. b. Nested distance estimate $D = 3$, $T = 3$, $\lambda = 10$.

**Figure 13 Lower and upper bounds and the nested distance between stochastic process and a tree.**

Differently from Figure 13, Figure 14a considers the scenario tree $\widetilde{\mathbb{P}}$ and demonstrates a *saturation effect* dependent on its branchiness factor $b$. In particular, the nested distance improvements are slowing down as the number of branches grows. In Figure 14a we observe that the nested distance decreases by half already after the first ten iterations. By this, the approximation of continuous-state stochastic processes by small scenario trees can, indeed, result in a fine point accuracy.

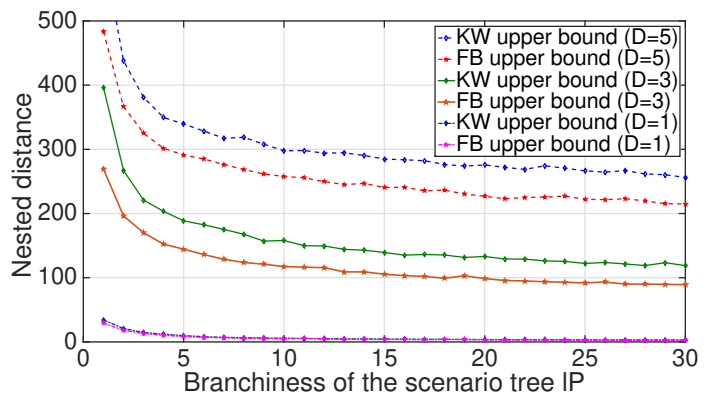

a. Branchiness.

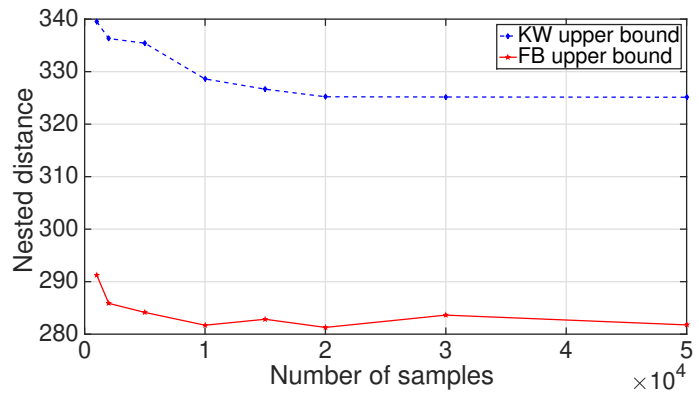

b. Number of samples.

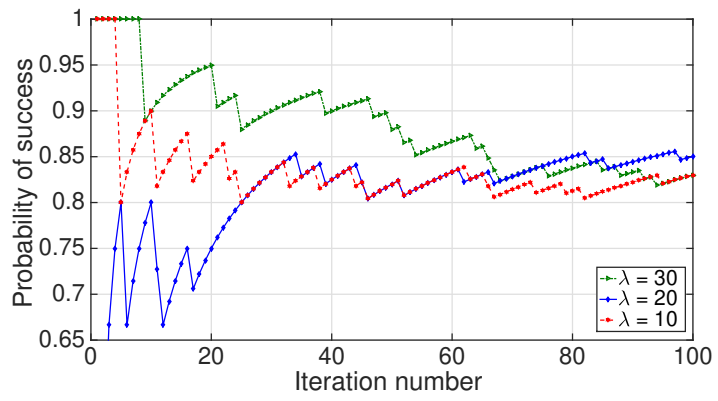

c. Success probability.

**Figure 14 Parameters influencing the nested distance and the probability to outperform the Kantorovich-Wasserstein quantization (i.e., success probability).**

While Figure 14b shows the convergence in the number of samples, Figure 14c demonstrates that the Algorithm 1 outperforms the stage-wise optimal quantization evaluated in line with the upper bound (7) in more than 85% of cases for the set of parameters described in Table 1. The underperformance in 15% of cases goes in line with the stochasticity of the methods and with the fact that Algorithm 1 is equivalent to stage-wise quantization if time stages are independent, i.e., if the covariance matrix is diagonal. If the covariance matrix is close to diagonal, we may also observe underperformance of Algorithm 1, whose *success probability* increases as the multiplier $\lambda$ grows (Figure 14c).

| | |
|---|---|
| Dimensionality, $D$ | $\mathbb{Z}[1,10]$ |
| Dependency multiplier, $\lambda$ | $\{10,20,30\}$ |
| Process mean, $\mu_t$ | $\lambda U[0,1]^D$, $\forall t$ |
| Process variance, $c_{t,t}$ | $\lambda U[0,1]$, $\forall t$ |
| Process covariance, $c_{s,t} = c_{t,s}$ | $\lambda U[0,1]$, $\forall s,t$ |
| Tree height, $T$ | $\mathbb{Z}[2,4]$ |
| Tree branchiness, $b$ | $\mathbb{Z}[2,5]$ |

**Table 1 Parameter values.**

Figures 9-14 demonstrate the results for the Gaussian distribution, which is $\|C_{t-1}^{-1}c^t\|$-Lipschitz continuous and for which all conditional distributions stay Gaussian in the form (17). Further, in Figure 15, we consider the case with the lognormal distribution $\log\mathcal{N}(\mu,C)$, which possesses similar properties as the Gaussian distribution and for which we can deduce the conditional form as $\exp(\xi_t|\xi^{t-1})$ with $\xi_t|\xi^{t-1} \sim \mathcal{N}\left(\mu_t + (\xi^{t-1} - \mu^{t-1})C_{t-1}^{-1}c^t, c_{t,t} - (c^t)^\top C_{t-1}^{-1}c^t\right)$ at each stage $t$. Here, instead of optimizing the quantizers directly for the lognormal distribution, we use an efficient approximation $\exp z_t^{(i)}, \forall i,t$ with the Gaussian quantizers $z_t^{(i)}, \forall i,t$ and their corresponding probabilities. By this, we analyze the drop in discretization quality assessing the opportunity to employ Gaussian quantizers for discretization of any probability distribution which can be represented as a combination of normal models. In Figure 15b, we observe a drop in the success probability (i.e., the probability that Algorithm 1 outperforms stage-wise optimal quantization) from ca. 85% in the Gaussian case to ca. 55% in the lognormal case. The reason for this is the approximation error in our efficient quantizer estimation approach. Though the conditional form of a lognormal distribution can be clearly deduced from the normal one, its optimal quantizers differ from $\exp z_t^{(i)}, \forall i,t$ and should be computed directly via Algorithm 1. Further, we observe an increase in the overall nested distance bound in Figure 15a in the lognormal case for the same tree structure.

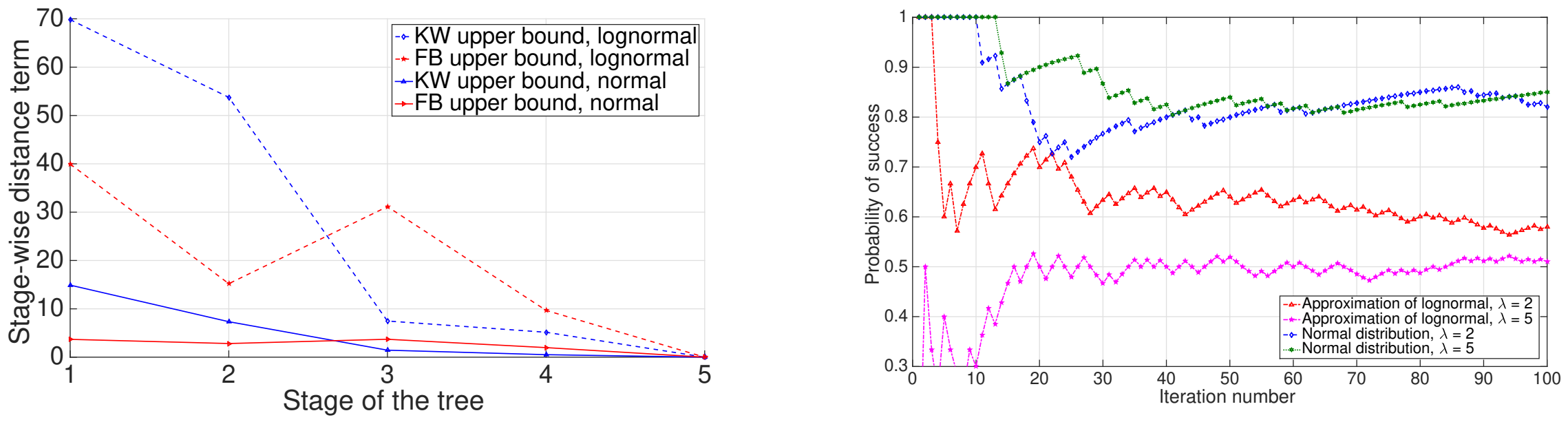

a. $D = 1, T = 5, \lambda = 1$. b. Probability comparison.

**Figure 15 Performance of the upper bound (7) for lognormal distribution.**

## 4. Multi-stage inventory control

To observe the performance of our algorithms on a real-world example, we solve of a multi-stage inventory control problem using the stage-wise optimal quantization and Algorithm 1 [8, 49, 63]. We consider a retailer facing time-varying demand uncertainties described via the stochastic process $\xi_1, ..., \xi_T$. The retailer places regular orders $x_{t-1}$ one period before the random demand realization is observed $\forall t = 1, ..., T$ and the cost for ordering one piece of a product is one. Unsold goods may be stored in the inventory with a storage loss $1 - l_t$. If the demand exceeds the inventory plus the newly arriving order, the demand has to be fulfilled by rapid orders (which are immediately delivered), for the price $h_t > 1$ per piece. The profit-maximization problem yields:

$$\begin{aligned} V = &\underset{x}{\text{maximize}} \quad \mathbb{E}\left[-\sum_{t=1}^{T}(x_{t-1} + h_t M_t) + l_T K_T\right] \\ &\text{subject to} \quad x_t \lhd \mathscr{F}_t,\ t = 1, ..., T, \\ &\qquad\qquad l_{t-1}K_{t-1} + x_{t-1} - \xi_t = K_t - M_t, \end{aligned} \tag{18}$$

where $K_t \geq 0$ is the inventory level right after all sales have been effectuated at time $t$ with $K_0 = 0$ and $M_t \geq 0$ is the corresponding shortage. The optimal profit is equal to $V + \mathbb{E}[\sum_{t=1}^{T} s_t \xi_t]$, where $s_t > 1$ denotes the selling price set by the retailer.

In a one-dimensional case, the closed-form solution can be derived for the problem (18) and, thus, the approximate solution can be compared with the true one [49, 63]. In a multi-dimensional case, the solution differs for independent and interdependent product demands. Indeed, for independent demands, the true optimal value is simply the sum of optimal values of one-dimensional problems due to the absence of a capacity constraint. For the multi-dimensional case with dependent demands, the explicit solution is unknown. The closed-form solution, i.e., the optimal value $V$ and the optimal decision $x_t$, $\forall t = 1, ..., T$ of the problem (18), can be written in the following form for the case with one product (see [49]):

$$\begin{aligned} V &= -\sum_{t=1}^{T} h_t \mathbb{E}(\xi_t) + \sum_{t=1}^{T}(h_t - 1)\mathbb{E}(\mathbb{AV@R}_{\beta_t}(\xi_t|\mathscr{F}_{t-1})) \\ x_t &= \mathbb{V@R}_{\beta_{t+1}}(\xi_{t+1}|\mathscr{F}_t) - l_t K_t, \end{aligned} \tag{19}$$

where $\beta_t = \frac{h_t - 1}{h_t - l_t}$, $K_t = \max\left\{\mathbb{V@R}_{\beta_t}(\xi_t|\mathscr{F}_{t-1}) - \xi_t,\, 0\right\}$, $\mathbb{V@R}_{\beta_t}(\xi) = F^{-1}(\beta_t)$, $\mathbb{AV@R}_{\beta_t}(\xi) = \frac{1}{\beta_t}\int_0^{\beta_t} F^{-1}(u)du$ and $F$ is the cumulative distribution function of a random variable $\xi$.

Further, we consider the Gaussian and the lognormal distributions to test the performance of Algorithm 1. In the Gaussian case, the optimal value can be written as

$$V = -\sum_{t=1}^{T} h_t\mu_t + \sum_{t=1}^{T}(h_t - 1)\left[\mu_t - \frac{1}{\beta_t}\sqrt{\frac{c_{t,t} - (c^t)^\top C_{t-1}^{-1} c^t}{2\pi}} \exp\left(-\frac{1}{2}\left[\Phi^{-1}(\beta_t)\right]^2\right)\right], \tag{20}$$

where $\Phi(\cdot)$ is the standard normal distribution function. Analogically, the lognormal demand implies

$$V = -\sum_{t=1}^{T} h_t \exp\left(\mu_t + \frac{c_{t,t}}{2}\right) + \sum_{t=1}^{T}(h_t - 1)\left[1 - \exp\left(\mu_t + \frac{c_{t,t}}{2}\right)\frac{\Phi\left(\Phi^{-1}(\beta_t) - \sqrt{c_{t,t} - (c^t)^\top C_{t-1}^{-1} c^t}\right)}{\beta_t}\right].$$

Both closed-form optimal values depend only on the mean vector and the covariance matrix of the random process $\xi = (\xi_1, ..., \xi_T)$ and play an important role for the comparison of different numerical algorithms. Approximating the stochastic process $\xi$ defined on the probability space $(\Omega, \mathscr{F}, P)$ by a scenario tree with the branchiness $b$, one can numerically calculate the solution of the inventory control problem using (i) a random (Monte-Carlo) scenario generation method, (ii) a stage-wise optimal quantization algorithm, as well as (iii) our forward-backward quantization method (Algorithm 1).

For the Gaussian and the lognormal cases, Figure 16 compares the optimal value $V$ with its numerical estimates obtained via the use of Algorithm 1 and the stage-wise optimal quantization. It also compares the convergence of these methods to the one for a random tree. As before, we increase the branchiness factors $b$ of the scenario tree and use Lemma 2 in the Appendix to observe the convergence. As one can see in Figure 16a, both methods converge to the optimal value with fine accuracy as the branchiness increases, while, as expected, the convergence of the random quantization method is much more volatile. The convergence is also observed in Figure 16b, where the demand follows the lognormal distribution.

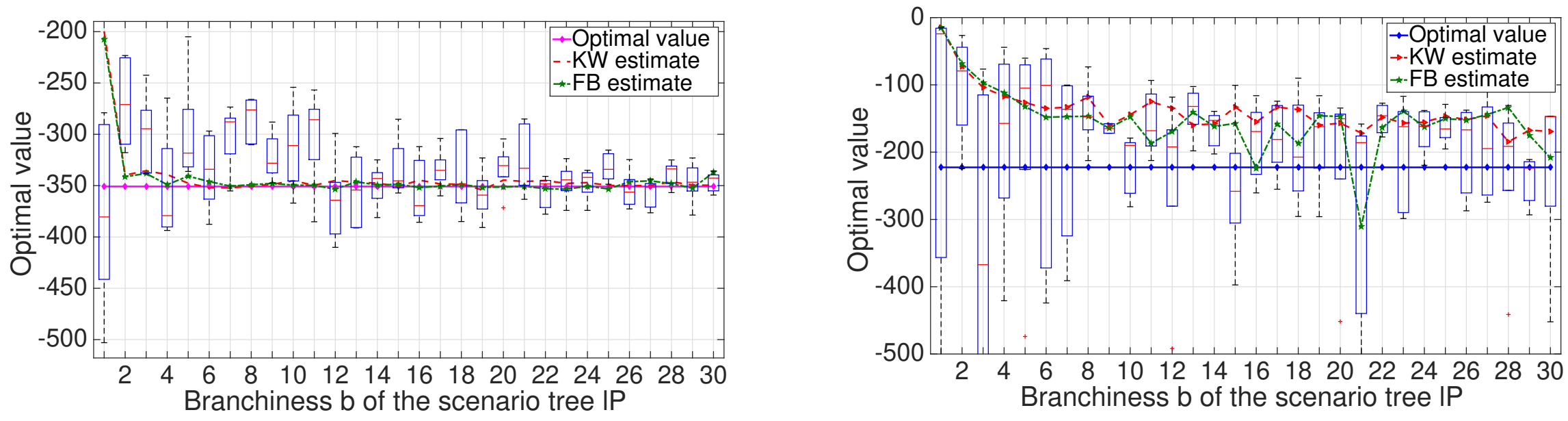

a. Gaussian distribution. b. Lognormal distribution.

**Figure 16 Convergence of the optimal value in one-dimensional inventory control.**

Further, the optimal values are linear in the number of time stages $T$ for both Gaussian and lognormal patterns because of the stationarity of demand processes in Figures 17a,b. As one can see in Figures 16b and 17b, the residual of the lognormal case is higher than that of the Gaussian one and is increasing in the number of stages $T$: this is expected due to the fact that we use a simple transformation of normal quantizers to lognormal ones by taking their exponents for both KW and FB methods (stage-wise optimal quantization and 1, respectively) and, thus, the result is the case described in Figure 15.

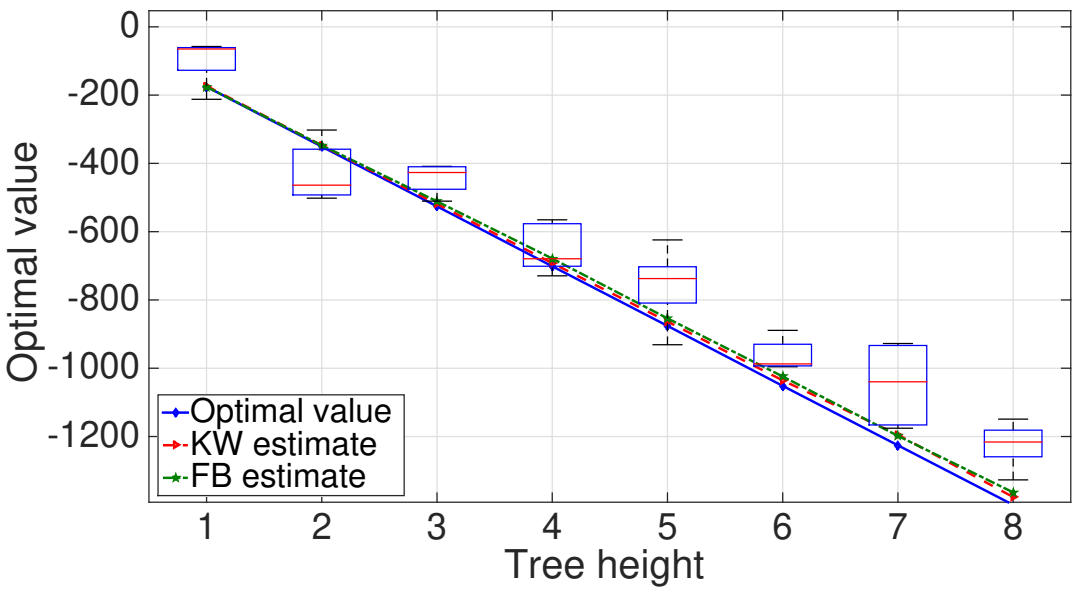

a. Gaussian distribution.

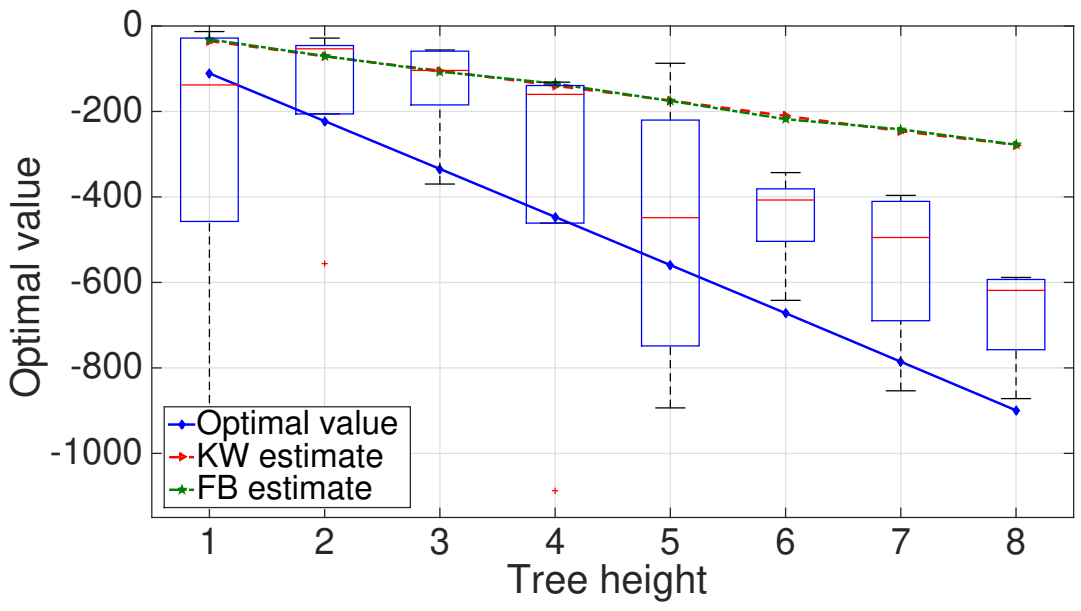

b. Lognormal distribution.

**Figure 17 Optimal value in one-dimensional inventory control with respect to the tree height.**

In Figure 18a,b, we analyze the average performance of (i) the Monte-Carlo tree sampling, (ii) the stage-wise optimal quantization and (iii) the forward-backward quantization in line with Algorithm 1 for the solution of the one-dimensional inventory control problem with normally distributed demand. For this, we consider stationary scenario processes with $\mu_s = \mu_t$, $\forall s,t \in \{1,...,T\}$ and a varying covariance matrix $C$ with $c_{t,t} = \mu_t$, which goes in line with demand approximations described in [61, 67]. We sample 50 inventory control problems with the parameters described in Table 2.

| | |
|---|---|
| Process mean, $\mu_t$ | $\mu_t = \mu_s = 100$, $\forall t,s$ |
| Process variance, $c_{t,t}$ | $c_{t,t} = \mu_t$ $\forall t$ |
| Process covariance, $c_{s,t} = c_{t,s}$ | $10U[0,1]$, $\forall s,t$ |
| Storage loss, $l_t$ | $0.1U[0,1]$, $\forall t$ |
| Rapid orders, $h_t$ | $1+U[0,1]$, $\forall t$ |
| Tree height, $T$ | $\mathbb{Z}[2,3]$ |
| Tree branchiness, $b$ | $\mathbb{Z}[1,30]$ |

**Table 2 Parameter values for the inventory control problem.**

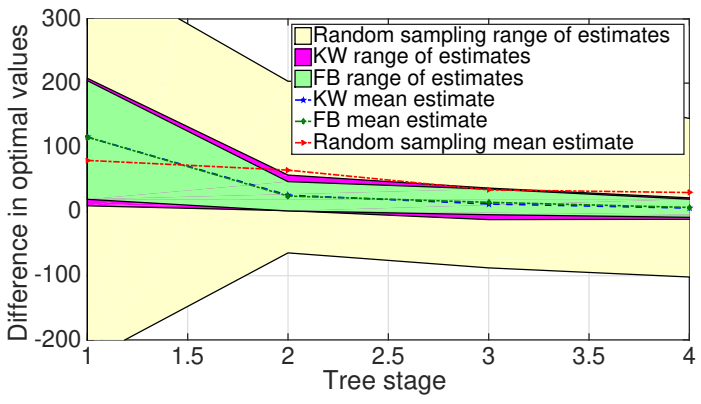

a. Range of estimates.

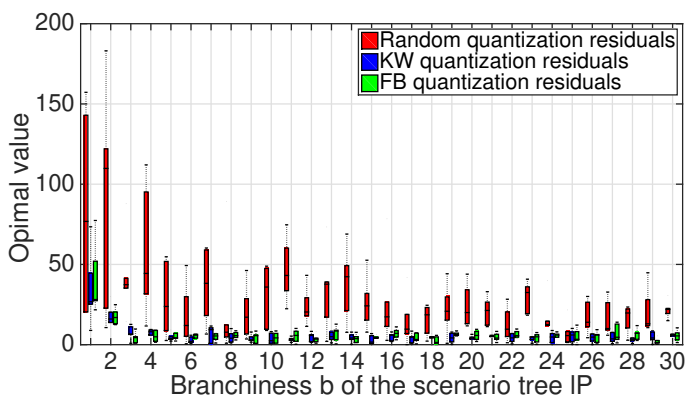

b. Residual value.

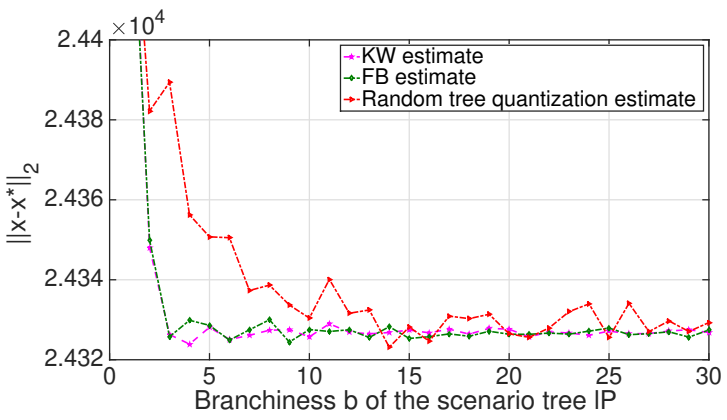

c. Residual in solutions.

**Figure 18 Average estimates in a one-dimensional inventory control problem.**

For small-size trees, the forward-backward quantization results in a narrower range of estimates for the inventory control problem than the stage-wise optimal quantization or, moreover, than the Monte-Carlo method. This implies higher approximation quality of Algorithm 1 (Figure 18a). For scenario trees with a larger branchiness factor, the performance of the stage-wise optimal quantization and Algorithm 1 start equalizing (Figure 18b,c): this goes in line with the fact that the discrete approximation obtained via Algorithm 1 approaches the continuous distribution $P_t$ if the number of branches in the scenario tree grows.

Next, we test different shapes of mean demand patterns, including constant, increasing, decreasing and a product life cycle curve (bell-shape). This is a step towards a natural extension of our work, which would account for possible *model uncertainties* (e.g., uncertainties in mean and variance of the demand distribution) but would require robust or distributionally robust techniques (see [35, 43]). In Figures 19 and 20, the optimal average order stays above the mean demand curve for every pattern, while the decision-maker mitigates peaks in demand by an increase in his inventory. Both considered quantization methods obtain reasonable estimates in a one-dimensional inventory control problem. The optimal value under the stationary demand (Figure 19a) is equal to $-1399.9$ in line with equation (20) for the Gaussian case. It is just 2.5% lower than the estimate obtained using Algorithm 1 (i.e., $-1365.4$) and ca. 3% lower than the estimate obtained via the stage-wise quantization (i.e., $-1358.5$). Next, the optimal value of the inventory control problem with the decreasing demand pattern (Figure 20b) is equal to $-6322.4$ and is ca. 1.9% lower than the estimate obtained using Algorithm 1 (i.e., $-6197.5$).

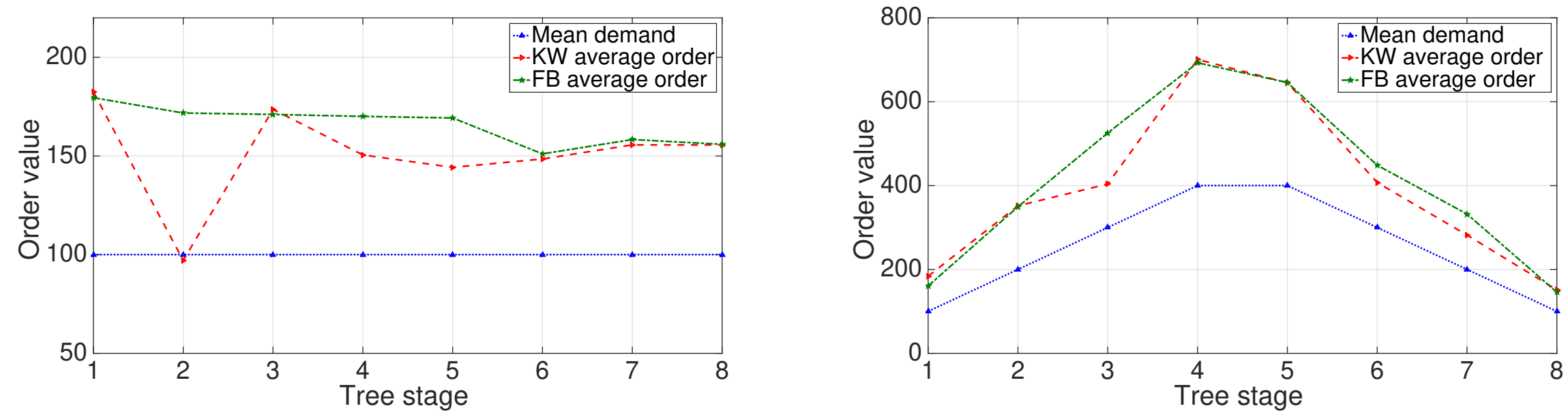

a. Stationary mean. b. Product life cycle.

**Figure 19 Ordering strategies of a unique product (1-dimensional case, constant or bell-shape demand).**

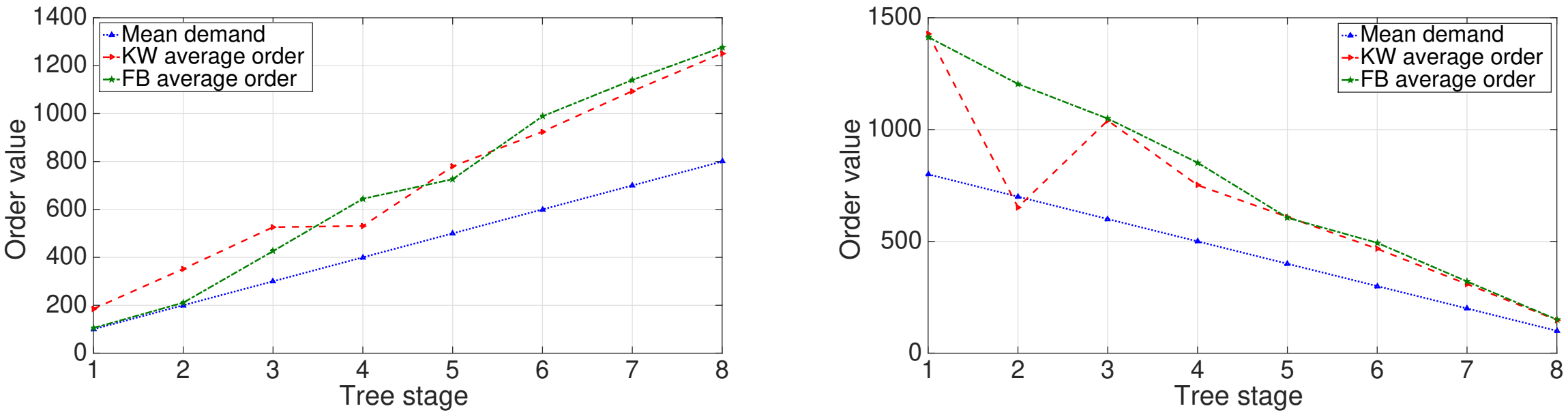

a. Increasing mean. b. Decreasing mean.

**Figure 20 Ordering strategies of a unique product (1-dimensional case, increasing or decreasing demand).**

As follows from equation (4), the optimal stock level of a unique product is defined by the $\beta_{t+1}$-quantile of the demand distribution and does not necessarily coincide with its mean (i.e., it coincides with the mean if $0.5h_t = 1 - 0.5l_t$ in the Gaussian case). This is observed in Figures 19 and 20 where the order quantity and, thus, the stock level are larger than the demand at each period. A similar result is observed for the multi-dimensional case with dependent products (Figure 21), for which the optimal order quantities exceed the mean demand.

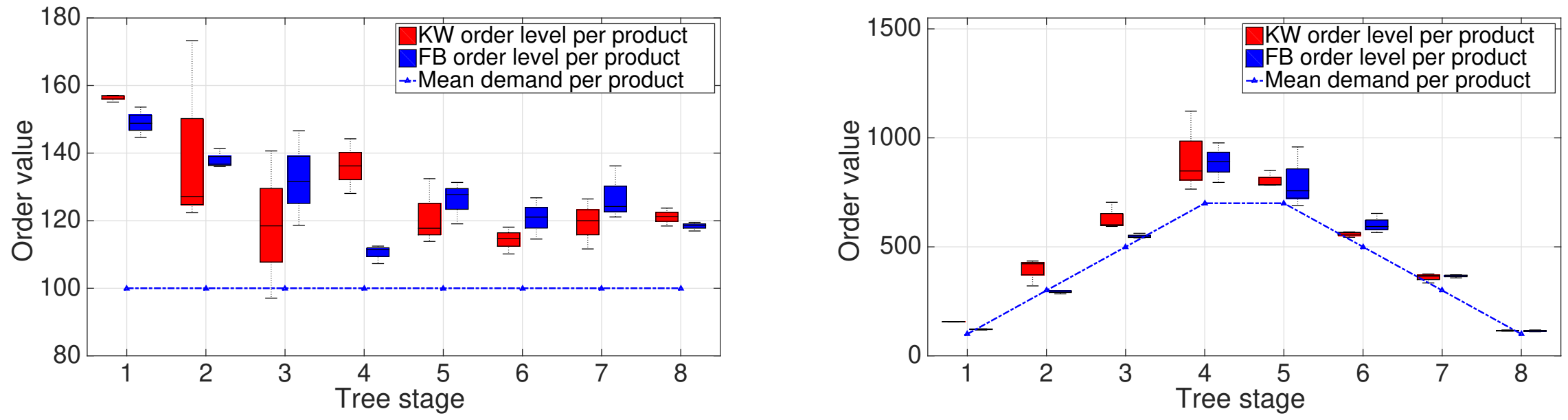

a. Optimal value. b. Order quantity.

**Figure 21 Ordering strategies in the presence of multiple interdependent products (3-dimensional case).**

Next, in Figure 22a,b, we demonstrate the dependence of the optimal value on model parameters for the inventory control problem with a unique product and a normally distributed demand. Modeling the urgent delivery price as $h_t = 1 + l_t$ and increasing the parameter $l_t$, we observe the convexity of the optimal value in $\frac{l_t}{1+l_t}$. Thus, even though the storage loss $1 - l_t$ decreases (i.e., $\frac{l_t}{1+l_t}$ increases), we observe a nadir in the retailer's profit due to an increase in the urgent delivery price (Figure 22a). A similar result is observed for the case with three interdependent products. Moreover, Figure 22c shows that an increase in the dependency multiplier $\lambda$ reduces the profit minimal value even further.

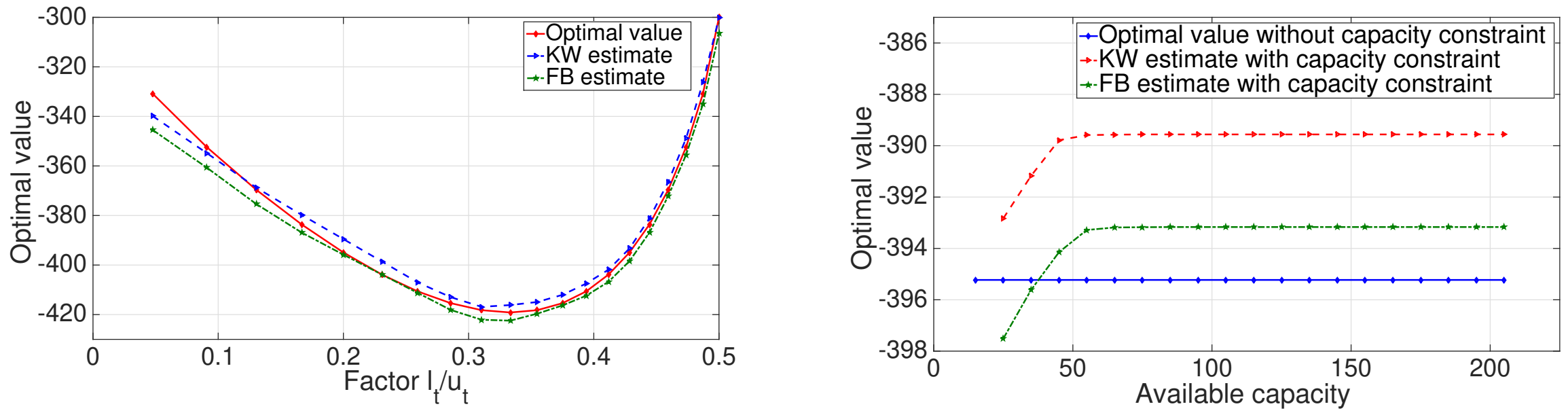

a. Storage loss v.s. urgent delivery, 1 product. b. Capacity constraint, 1 product.

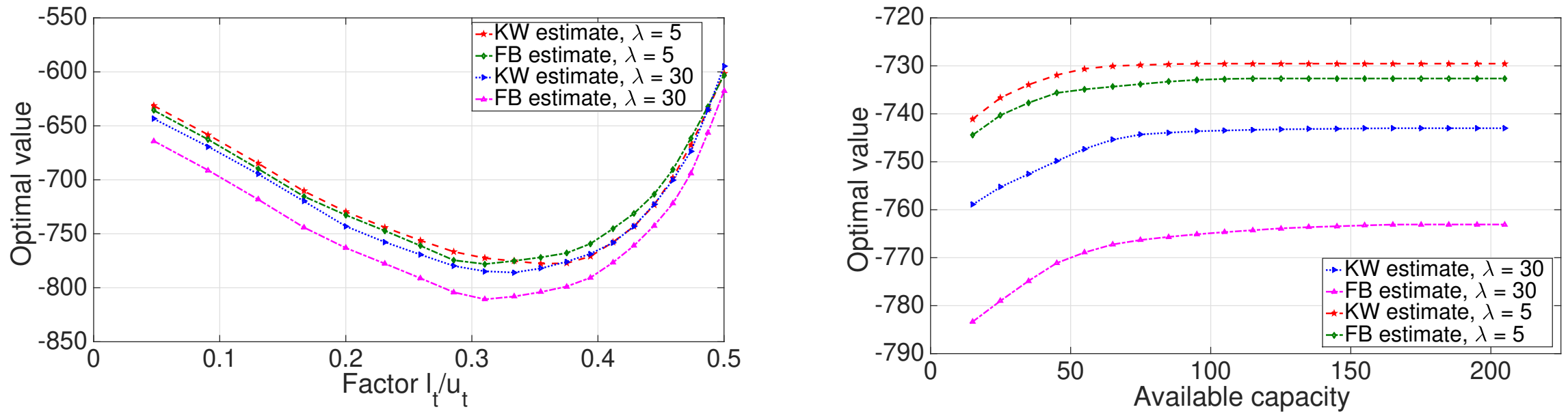

c. Storage loss v.s. urgent delivery, 3 products. d. Capacity constraint, 3 products.

**Figure 22 Optimal value of the inventory control problem dependent on the model parameters.**

Including the capacity constraint $l_{t-1}K_{t-1} + x_t < \hat{C}$ to the inventory control problem and decreasing the capacity $\hat{C}$ obviously results in a drop of the optimal value (Figure 22b,c). Differently, increasing the capacity at the retailer, we demonstrate the convergence of the optimal value. Clearly, for small-size scenario trees, the convergence attains the value in some distance to the true solution of the continuous inventory control problem, e.g., Figure 22b,c corresponds to the scenario tree with branchiness $b = 5$ and height $T = 2$.

## 5. Conclusion

In this article, we study numerical methods for the solution of multi-stage stochastic optimization problems by the use of scenario approximation techniques. Working in a purely distributional setup, we approximate continuous-state stochastic processes by finitely valued scenario trees, developing a forward-backward multi-stage scenario generation algorithm. Measuring the approximation error via the use of nested distances, we introduce and minimize an upper bound effectively combining information about the past with available information about feasible future scenarios. Compared to other existing bounds, our upper bound does not neglect the fact that the set of feasible future scenarios at any particular tree node is smaller than the set of all future scenarios. Also, introducing a concept of clairvoyant trees, we prove that their particular sequence improves the existing lower bound on the nested distance.

Introducing probabilistic information about feasible future subtrees explicitly to the scenario generation, we make a step towards new methods for distributional quantization specifically adjusted to multi-stage optimization. We use the proposed upper bound and develop an algorithm for its minimization. We numerically demonstrate that the resulting quantizers provide better guarantees for the minimal nested distance than the quantizers obtained by the stage-wise minimization of the Kantorovich-Wasserstein distance or, moreover, than the random (Monte-Carlo) sampling. Furthermore, we show that even small-size trees can result in a fine approximation of stochastic processes given that the quantization takes the interdependencies between time stages explicitly into account. In particular, this is important for distributionally robust cases, where the center of the ambiguity set must be a small-size but also a fine scenario tree.

The developed scenario approximation method can be applied in a huge variety of areas: starting with financial planning and inventory control, possible applications include topics of energy production and trading, electricity generation planning, pension fund management, supply chain management and similar fields. Importantly, in optimization problems of asset management, one needs to impose no-arbitrage constraints in the optimization problem for scenario generation. This is necessary for consistency with financial asset pricing theory (see [27, 32]). In our future research, we would like to focus on forward-backward splitting techniques for stage-wise problems in optimal transport and on scenario generation methods for ambiguous distributions, i.e., for distributions which vary not only in time but also in some uncertainty set (see [35]).

## 6. Acknowledgements

We are very grateful to the Austrian Science Fund (FWF Der Wissenschaftsfonds), grant number J3674-N26 for funding our research on accuracy and efficiency of approximation methods in multi-stage stochastic optimization. Furthermore, we would like to thank Prof. Dr. Daniel Kuhn for very valuable remarks about distributional properties of scenario trees.

# Appendix

EXAMPLE 2. Consider nested distributions $\mathbb{P}$ and $\widetilde{\mathbb{P}}$ and suppose that one of them becomes clairvoyant at stage $T-1$, i.e., we obtain $\mathbb{P}_{c(T-1)}$. Next, we would like to test if the following inequality holds:

$$dl(\mathbb{P}_{c(T-1)}, \widetilde{\mathbb{P}}) \leq dl(\mathbb{P}, \widetilde{\mathbb{P}}), \tag{21}$$

that is similar to the lower bound (8), but with just one clairvoyant tree. Clearly, if nested distributions $\mathbb{P}$ and $\widetilde{\mathbb{P}}$ are equal ($\mathbb{P} = \widetilde{\mathbb{P}}$), the nested distance between them is equal to zero, i.e., $dl(\mathbb{P}, \mathbb{P}) = 0$. Thus, if one of the nested distributions becomes clairvoyant at stage $T-1$, the distance $dl(\mathbb{P}_{c(T-1)}, \mathbb{P})$ should also be equal to zero in order to satisfy (21). Consider the following tree:

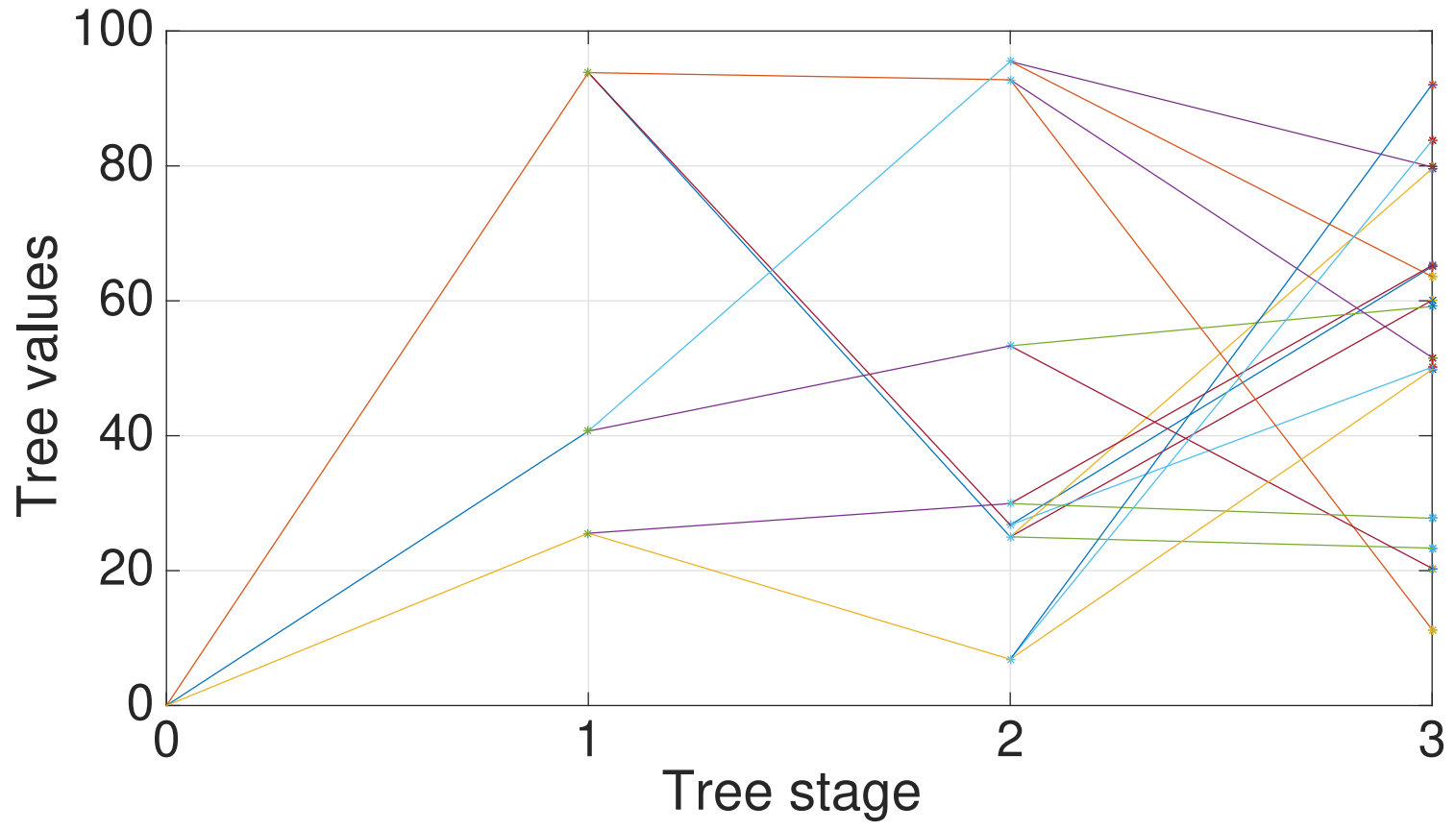

**Figure 23 Tree structure with the following probability vector $p$ and the corresponding values $v$:**
$p$ = [1 0.2 0.3 0.5 0.5 0.5 0.6 0.2 0.2 0.7 0.3 0.5 0.5 0.1 0.9 0.2 0.8 0.1 0.2 0.7 0.5 0.5 0.2 0.3 0.5 0.9 0.1]
$v$ = [0 40.68 93.83 25.54 53.32 95.48 26.77 25.01 92.77 6.86 29.94 59.16 20.33 63.59 79.84 50.17 65.08 79.6 23.34 60.08 11.25 51.58 83.78 92.08 49.82 27.76 65.25].

Changing the structure of the tree so that it becomes clairvoyant at stage $T-1$ and calculating necessary nested distances, we obtain the counterexample for the inequality (21), i.e., $11.3919 = dl(\mathbb{P}_{c(T-1)}, \mathbb{P}) > dl(\mathbb{P}, \mathbb{P}) = 0$ for the given tree. Thus, a nested distance $dl(\mathbb{P}_{c(T-1)}, \widetilde{\mathbb{P}})$ is not a lower bound for the distance $dl(\mathbb{P}, \widetilde{\mathbb{P}})$. The distance $dl(\mathbb{P}_{c(T-1)}, \widetilde{\mathbb{P}})$ is also not an upper bound for the nested distance $dl(\mathbb{P}, \widetilde{\mathbb{P}})$. To see this, let $\widetilde{\mathbb{P}} = \mathbb{P}_{c(T-1)}$. The distance $dl(\mathbb{P}, \mathbb{P}_{c(T-1)})$ is equal to 11.3919, as shown above. Nevertheless, the nested distance $dl(\mathbb{P}_{c(T-1)}, \mathbb{P}_{c(T-1)})$ is obviously equal to zero, i.e., $0 = dl(\mathbb{P}_{c(T-1)}, \mathbb{P}_{c(T-1)}) < dl(\mathbb{P}, \mathbb{P}_{c(T-1)})$, where we again change just one tree structure to the clairvoyant one. □

THEOREM 2. (see [20, 50, 71]) Suppose that the $D$-dimensional distribution $P$ has a density $f$ with $\int |u|^{1+\delta} f(u)du < \infty$ for some $\delta > 0$ and let

$$\left[ \int f(x)^{\frac{D}{D+1}} dx \right]^{\frac{D+1}{D}} \leq c,$$

where $c$ is a constant. Then, $\exists \widetilde{P}^N$, such that $d(P, \widetilde{P}^N) \leq cN^{-\frac{1}{D}}$, where we denote by $\widetilde{P}^N$ a discrete approximation of $P$ with $N$ values. Furthermore,

$$d(P, \widetilde{P}^N) \to 0, \text{ if } N \to \infty.$$

**Proof:** The result of the theorem follows from the Zador-Gersho formula (see [20, 50, 71]). □

LEMMA 2. (Limit principle) Let the nested distribution $\mathbb{P}$ be approximated by a scenario tree $\widetilde{\mathbb{P}}$ with $b_t$ number of branches outgoing from every node at stages $t = 1,...,T-1$ of the tree. Let the $(K_2,...,K_T)$-Lipschitz property hold for the multivariate distribution $P^{1:T}$, that is dissected into the chain $P_1$, $P_2$, ..., $P_T$ of conditional distributions. If *branchiness factors* $b_1,...,b_T$ are monotonically increasing and the conditions of Theorem 2 are satisfied for each $P_t$, then the nested distance $dl(\mathbb{P},\widetilde{\mathbb{P}})$ is converging to zero.

**Proof:** Denote by $\widetilde{P}_t$ the discrete approximation of $P_t$ sitting on $b_t$ points. As the conditions of Theorem 2 are satisfied, $d_{KA}(P_t,\widetilde{P}_t^{b_t}) \to 0$ with the branchiness factor $b_t \to \infty\ \forall t = 1,...,T$. Further, as the $(K_2,...,K_T)$-Lipschitz property holds for the distribution $P_t$, we can write

$$dl(\mathbb{P},\widetilde{\mathbb{P}}) \leq \sum_{t=1}^{T} d_{KA}(P_t,\widetilde{P}_t^{b_t}) \prod_{s=t+1}^{T} (K_s+1).$$

This implies that the upper bound for the nested distance is converging to zero when the branchiness of the tree is increasing and, hence, $dl(\mathbb{P},\widetilde{\mathbb{P}}) \to 0$ as $b_t \to \infty$, $\forall t$. □

LEMMA 3. (Nested Distance Convergence) Let $\mathbb{P}$ be a nested distribution with continuous multivariate distribution $P^{1:T}$ dissected into the chain $P_1$, $P_2$, ..., $P_T$ of conditional distributions and suppose that there is another nested distribution $\widetilde{\mathbb{P}}$ with a finite tree structure. If conditions of Lemma 2 are satisfied, then there exists such finite nested distribution $\mathbb{Q}_b$ sitting on a tree with branchiness factor $b\ \forall t = 1,...,T$, for which the following convergence result holds: $\exists \mathbb{Q}_b:\ \lim_{b\to\infty} dl(\mathbb{P},\mathbb{Q}_b) = 0,\ \lim_{b\to\infty} dl(\mathbb{Q}_b,\widetilde{\mathbb{P}}) = dl(\mathbb{P},\widetilde{\mathbb{P}})$.

**Proof:** The result directly follows from the triangle inequality $dl(\mathbb{P},\widetilde{\mathbb{P}}) \leq dl(\mathbb{P},\mathbb{Q}_b) + dl(\mathbb{Q}_b,\widetilde{\mathbb{P}})$ and Lemma 2, under which the nested distance $dl(\mathbb{P},\mathbb{Q}_b)$ converges to zero when the branchiness factor $b$ approaches infinity: $b \to \infty \Rightarrow dl(\mathbb{P},\mathbb{Q}_b) \to 0 \Rightarrow dl(\mathbb{Q}_b,\widetilde{\mathbb{P}}) \to dl(\mathbb{P},\widetilde{\mathbb{P}})$. □

## 6.1. Stochastic gradient algorithms

Let $f(z,\xi)$ be a function, where $\xi \in \mathbb{R}^D$ is a random variable and $z \in \mathbb{R}^{D'}$ is a decision variable. We assume that $F(z) := \mathbb{E}[f(z,\xi)]$ exists and is differentiable in $z$. Moreover, we assume that there is a *stochastic gradient* $g(z,\xi)$ such that $\mathbb{E}[g(z,\xi)] = \nabla_z F(z)$. For finding a (local) minimum of $F$, the *stochastic gradient algorithm* is based on a sequence of independent copies $\xi_1,\xi_2,\dots$ of $\xi$ and approximates the minimizer of $F$ by the iteration $z_{k+1} = z_k - t_k g(z_k,\xi_k)$, where $t_k$ is a sequence of step-sizes, e.g., the diminishing steps such that $\sum_k t_k = \infty$ and $\sum_n t_k^2 < \infty$. Convergence of this algorithm was studied by numerous authors in the years 1950-1980, see e.g., [42, 54, 58]. If the problem is to find $\min\{F(z) : z \in \Xi\}$ with a convex constraint set $\Xi$, then the method is modified to the projected gradient algorithm:

$$z_{k+1} = P_\Xi(z_k - t_k g(z_k,\xi_k))$$

where $P_\Xi$ is the projection onto the convex constraint set $\Xi$. A special case (and this is the case we consider in this paper) is given if the constraint set is a linear space, i.e., if

$$\Xi = \{z \in \mathbb{R}^{D'} : Az = 0\},$$

where $A$ is a $[D'' \times D']$ matrix with $D'' < D'$ and $\text{rank}(A) = D''$. The specific structure of the constraint set allows us to compute the projection operator using the orthogonal projector. Specifically, the projection on $\mathbb{X}$ is given by the symmetric matrix $P_A = I_{D'} - A^T(AA^T)^{-1}A$, where $I_{D'}$ is the identity matrix of the size $D' \times D'$. Assuming that $z_k \in \Xi$, the stochastic gradient algorithm is

$$z_{k+1} = P_A(z_k - t_k g(z_k,\xi_k)) = z_k - t_k P_A(g(z_k,\xi_k))$$

ensuring that also $z_{k+1} \in \Xi$.